%% file: multitm-preprint.tex
\newcount\Comments  
\documentclass[letterpaper,10pt]{scrartcl}

\usepackage[utf8]{inputenc}
\usepackage[english]{babel}
\usepackage[automark]{scrpage2}
\usepackage{amsfonts}
\usepackage{amsmath}
\usepackage{amssymb}
\usepackage{amsthm}
\usepackage{color}
\usepackage{graphicx}
\usepackage{sidecap}
\usepackage{multirow}
\usepackage{booktabs}
\usepackage{bm}
\usepackage{fullpage}

\usepackage{tabularx}
\usepackage[font=small]{caption}
\usepackage{subfig}

\usepackage{longtable}
\usepackage{rotating}

\usepackage{algorithm}
\usepackage{algpseudocode}
\usepackage{lineno,hyperref}

\usepackage{tabularx}
\usepackage[font=small]{caption}
\usepackage{subfig}

\usepackage{longtable}
\usepackage{rotating}

\usepackage[textsize=tiny]{todonotes}

\usepackage{algorithm}
\usepackage{algpseudocode}

\graphicspath{{figures/}}
\DeclareGraphicsExtensions{.pdf,.eps,.png,.jpg,.jpeg}

\newcommand{\bfx}{\boldsymbol x}

\newcommand{\bfbeta}{\boldsymbol \beta}

\newcommand{\bfy}{\boldsymbol y}

\newcommand{\Jcal}{\mathcal{J}}

\newcommand{\Pcal}{\mathcal{P}}

\newcommand{\Tcal}{\mathcal{T}}

\newcommand{\Ycal}{\mathcal{Y}}

\newcommand{\bfepsilon}{\boldsymbol \epsilon}
\newcommand{\bfSigma}{\boldsymbol \Sigma}
\newcommand{\bfj}{\boldsymbol j}

\newcommand{\Var}{\operatorname{Var}}

\newcommand{\kibitz}[2]{\ifnum\Comments=1\textcolor{#1}{#2}\fi}

\newcommand{\refP}{\eta}
\newcommand{\pb}{\widetilde{\refP}}
\newcommand{\pbLow}{\widehat{\refP}}
\newcommand{\post}{\pi}
\newcommand{\postLow}{\widehat{\post}}
\newcommand{\mapT}{T}
\newcommand{\mapTApprox}{\widetilde{\mapT}}
\newcommand{\mapTApproxLow}{\widehat{\mapT}}
\newcommand{\mapTInv}{S}
\newcommand{\mapTInvApproxLow}{\widehat{\mapTInv}}
\newcommand{\sampPi}{\theta}
\newcommand{\bfSampPi}{\boldsymbol{\sampPi}}
\newcommand{\sampR}{\vartheta}
\newcommand{\bfSampR}{\boldsymbol{\sampR}}

\ihead[]{}
\ohead[]{}
\ifoot[]{}
\ofoot[]{}

\newenvironment{keywords}%
   {\begin{trivlist}\item[]{\bfseries\sffamily Keywords:}\ }
   {\end{trivlist}}

\ifpdf
\hypersetup{
  pdftitle={A transport-based multifidelity preconditioner for Markov chain Monte Carlo},
  pdfauthor={Benjamin Peherstorfer, Youssef Marzouk}
}
\fi

\title{A transport-based multifidelity preconditioner for Markov chain Monte Carlo\thanks{The second author acknowledges support of the AFOSR MURI on multi-information sources of multi-physics systems under Award Number FA9550-15-1-0038.}
}

\author{Benjamin Peherstorfer\footnote{Courant Institute of Mathematical Sciences, New York University, New York, NY 10012 (pehersto@cims.nyu.edu)} \and
        Youssef Marzouk\footnote{Massachusetts Institute of Technology, Cambridge, MA 02139 (ymarz@mit.edu)} 
}

\date{August 2018}

\begin{document}

\maketitle

\begin{abstract}
Markov chain Monte Carlo (MCMC) sampling of posterior distributions arising in Bayesian inverse problems is challenging when evaluations of the forward model are computationally expensive. Replacing the forward model with a low-cost, low-fidelity model often significantly reduces computational cost; however, employing a low-fidelity model alone means that the stationary distribution of the MCMC chain is the posterior distribution corresponding to the low-fidelity model, rather than the original posterior distribution corresponding to the high-fidelity model. We propose a multifidelity approach that combines, rather than replaces, the high-fidelity model with a low-fidelity model. First, the low-fidelity model is used to construct a transport map that deterministically couples a reference Gaussian distribution with an approximation of the low-fidelity posterior. Then, the high-fidelity posterior distribution is explored
using a non-Gaussian proposal distribution derived from the transport map. This multifidelity ``preconditioned'' MCMC approach seeks efficient sampling via a proposal that is explicitly tailored to the posterior at hand and that is constructed efficiently with the low-fidelity model. By relying on the low-fidelity model only to construct the proposal distribution, our approach guarantees that the stationary distribution of the MCMC chain is the high-fidelity posterior. In our numerical examples, our multifidelity approach achieves significant speedups compared to single-fidelity MCMC sampling methods.

\end{abstract}


\section{Introduction}
Bayesian inference provides a framework to quantify uncertainties in the solutions of inverse problems \cite{KaipoBook,ANU:7701764,TarantolaBook}. The Bayesian approach to inverse problems combines observed data, a forward model that maps parameters to observations, a prior distribution on the parameters, and a statistical model for the mismatch between model predictions and observations to define the \emph{posterior distribution} via Bayes' theorem. The posterior distribution characterizes the parameter values and their uncertainties, given these  ingredients. Practically ``solving'' a Bayesian inverse problem, however, entails exploring the posterior distribution, e.g., computing posterior expectations. A flexible and widely used approach for exploring posterior distributions is to draw samples with Markov chain Monte Carlo (MCMC) methods \cite{tierney1994,Gilks1996}. Using these methods, the forward model typically must be evaluated multiple times at different parameter values for each sample that is drawn, such that MCMC sampling quickly becomes computationally infeasible if each forward model solve is expensive. 

In this paper, we propose a multifidelity preconditioner to increase the efficiency of MCMC sampling. Our multifidelity approach exploits low-cost, low-fidelity models to construct a proposal distribution that approximates the posterior distribution at hand, then uses this proposal distribution to perform MCMC sampling of the original (high-fidelity) posterior distribution. The proposal distribution is derived from a \emph{transport map} that transforms the potentially complex posterior distribution into another distribution from which samples can be drawn more easily; in particular, we seek a map that transforms the posterior into a more Gaussian and more isotropic distribution. See \cite{ELMOSELHY20127815,Marzouk2016,ParnoTMap} for an introduction to transport maps in the context of Bayesian inverse problems. The key idea in the present work is to construct an invertible map using low-fidelity models, but to apply it to the high-fidelity posterior. The map then acts as a \emph{preconditioner} for the high-fidelity posterior, preserving information (due to its invertibility) while enabling MCMC sampling to proceed more efficiently.
An alternative but equivalent perspective is that pushing a simple proposal through the inverse of this transport map---for instance, using a Metropolis independence sampler with a standard Gaussian proposal---yields potentially non-Gaussian and tailored proposals that can be efficient for the high-fidelity posterior. The end result is the same: more efficient MCMC sampling, allowing the number of high-fidelity forward model evaluations to be reduced.
Since the low-fidelity model is used only for preconditioning (or equivalently, constructing the proposal), the stationary distribution of the chain obtained with our multifidelity approach is guaranteed to be the posterior distribution corresponding to the high-fidelity model.

There is a long tradition of exploiting low-fidelity models to speed up MCMC sampling for Bayesian inverse problems. There is work \cite{WANG200515,WRCR:WRCR9576,MARZOUK2007560,MARZOUK20091862,Lieberman2010,CHEN201584} that \emph{replaces} the high-fidelity model with a low-fidelity model (see also the survey \cite{Frangos}); because of this replacement, however, samples are drawn from a different distribution---the posterior induced by the low-fidelity model. Under certain assumptions, it has been shown that the posterior distribution induced by the low-fidelity model converges to the original/high-fidelity posterior distribution as the low-fidelity approximation is refined \cite{MarzoukXiu2009,CotterDashtiStuartSINUM2010}. Other approaches \emph{adapt} low-fidelity models \cite{NME:NME4748,YoussefJinglai} over a finite interval of posterior exploration, or quantify the error introduced by sampling the low-fidelity posterior distribution \cite{ROMES,ManzoniInverse}. Yet another family of approaches incrementally and infinitely refines approximations of the forward model on-the-fly during MCMC sampling \cite{ConradLocal,ConradDavisetalJUQ2018}; under appropriate conditions, these schemes guarantee that the MCMC chain asymptotically samples the high-fidelity posterior distribution. 

Instead of replacing the high-fidelity model with low-fidelity models, multifidelity methods combine high- and low-fidelity models. The aim is to leverage low-fidelity models for speeding up computations while allowing occasional recourse to the high-fidelity model to establish accuracy guarantees \cite{PWG17MultiSurvey}. A variety of multifidelity methods have been developed for uncertainty propagation \cite{boyaval_variance_2010,boyaval_fast_2012,NME:NME4761,Eldred2016,CHEN2013233,Peherstorfer16MFIS,Peherstorfer15Multi,MFGSA,doi:10.1137/17M1122992}; see the survey \cite{PWG17MultiSurvey}. For the solution of Bayesian inverse problems, there are multi-stage MCMC methods that aim to reduce the number of high-fidelity model evaluations by first screening proposed moves with low-fidelity models  \cite{MultiStageMCMC,Efendiev}. Another line of work builds on hierarchies of low-fidelity models, typically derived from different discretizations of partial differential equations (PDEs) underlying the high-fidelity model, to reduce sampling costs \cite{BESKOS20171417,HierMCMC,2017arXiv170909763L}. 

An alternative to multi-stage and hierarchical methods for increasing the efficiency of MCMC sampling is to use notions of transport to construct more effective proposal distributions. Effective proposals in MCMC should reflect the local or global geometry of the target distribution \cite{Gilks1996,MCMCHandbook}. In keeping with this idea, Parno et al.~\cite{ParnoTMap} use transport maps to precondition MCMC sampling. As described earlier, preconditioning involves constructing transport maps that ``Gaussianize'' the target distribution, such that it can be sampled more effectively by standard MCMC algorithms; these maps thus encode the geometry of the target. 
Parno et al.~\cite{ParnoTMap} build and refine such transport maps in an online fashion during MCMC sampling, as more and more samples are obtained. This approach can be seen as a form of \emph{adaptive MCMC} \cite{haario2001,doi:10.1198/jcgs.2009.06134}, wherein a (non-Gaussian) proposal distribution is adapted as the MCMC sampling proceeds. The transport maps in \cite{ParnoTMap} are constructed via the solution of a convex and separable optimization problem, which is simple and fast to obtain numerically. However, this sampling approach faces some of the usual pitfalls of adaptive MCMC. One of these issues is that a certain amount of initial mixing is necessary for the adaptation to be effective, because the online-adapted transport map depends on past samples. 

The transport approach of Moselhy et al.~\cite{ELMOSELHY20127815} instead follows an offline/online decomposition of the computation. In an offline phase, a reference distribution is selected from which independent samples can be drawn cheaply and then a transport map is constructed that pushes forward the reference distribution to the posterior distribution. Then, in the online phase, the transport map is used to transform samples from the reference distribution into samples from the posterior distribution. If a large number of samples are transformed in the online phase, then the one-time high cost of constructing the transport map offline is compensated. The optimization problem for constructing the transport map can employ derivative information from the posterior density, and thus yields accurate maps even if the posterior is concentrated. Furthermore, no samples of the posterior distribution are necessary to solve the optimization problem, in contrast with adaptive MCMC techniques such as \cite{ParnoTMap}. However, evaluations of the objective of the optimization problem entail evaluations of the unnormalized posterior density and thus of the forward model; moreover, this problem is in general not separable across dimensions and not convex (unless the posterior density is log-concave). As a result, map construction in \cite{ELMOSELHY20127815} is typically far more computationally demanding than  map construction in \cite{ParnoTMap}.
Furthermore, since the construction of the transport map involves numerical approximations, the pushforward of the reference distribution by the map in general only \emph{approximates} the posterior distribution. 

We propose a multifidelity approach that combines several of the advantages of the approaches introduced by Parno \cite{ParnoTMap} and Moselhy \cite{ELMOSELHY20127815}. We construct transport maps building on the offline/online approach of \cite{ELMOSELHY20127815} and mitigate the high computational costs of solving the optimization problem by relying on low-cost, low-fidelity models. Our approach can exploit a wide range of low-fidelity models, including projection-based reduced models \cite{SirovichMethodOfSnapshots,RozzaPateraSurvey,SerkanReview,BennerSIREV,ChenSurvey}, data-fit interpolation and regression models \cite{forrester_recent_2009}, machine-learning-based models \cite{SVM,VKOGA,2018arXiv180709575S,VapnikBook}, and simplified-physics models \cite{allaire_mathematical_2014,ng_monte_2015}. Then, in the online phase, the transport map is used to precondition MCMC sampling of the high-fidelity posterior, as in \cite{ParnoTMap}. The corresponding MCMC scheme is ergodic for the high-fidelity posterior. In other words, the Metropolis step corrects errors that otherwise would be introduced by relying on the transport map alone to push forward the reference to the posterior.
Thus, we obtain a multifidelity approach that uses low-fidelity models to speed up computations while making occasional recourse to the high-fidelity model to establish convergence, in the sense that the stationary distribution of our MCMC chain is the high-fidelity posterior distribution. In contrast to \cite{ParnoTMap}, our MCMC algorithm is not adaptive, because the map is built once offline and then stays fixed during MCMC sampling. In particular, no samples of the high-fidelity posterior distribution are needed to construct the transport maps.

Section~\ref{sec:Prelim} describes the problem setup and briefly reviews transport maps in the context of Bayesian inverse problems. In Section~\ref{sec:MFTMap}, we introduce our multifidelity preconditioner and discuss it in the context of sampling with the Metropolis-Hastings algorithm. Section~\ref{sec:NumExp} demonstrates our multifidelity approach on two examples where we achieve significant speedups compared to using the high-fidelity model alone. Conclusions are given in Section~\ref{sec:Conc}.

\section{Preliminaries}
\label{sec:Prelim}
Section~\ref{sec:Prelim:BIP} and Section~\ref{sec:Prelim:MH} define our Bayesian inverse problem setting and describe the Metropolis-Hastings algorithm for MCMC sampling. We refer to, e.g., \cite{KaipoBook,ANU:7701764,TarantolaBook}, for details on Bayesian approaches to inverse problems. Section~\ref{sec:Prelim:TMap} discusses transport maps for coupling probability distributions in the context of Bayesian inverse problems. The problem description is given in Section~\ref{sec:Prelim:Problem}.

\subsection{Bayesian inverse problems}
\label{sec:Prelim:BIP}
Consider the high-fidelity forward model (parameter-to-observable map) $G: \Pcal \to \Ycal$, with parameter $\bfSampPi \in \Pcal$, where $\Pcal \subseteq \mathbb{R}^d$, and observable $\bfy \in \Ycal$, where $\Ycal \subseteq \mathbb{R}^{d^{\prime}}$. Thus, the parameter $\bfSampPi$ and the observable $\bfy$ are $d$-dimensional and $d^{\prime}$-dimensional vectors, respectively. Consider now observed data
\[
\bfy = G(\bfSampPi^*) + \bfepsilon\,,
\]
corresponding to some parameter value $\bfSampPi^* \in \Pcal$. The noise $\bfepsilon$ is assumed to be a realization of a zero-mean Gaussian random variable with covariance matrix $\bfSigma_{\epsilon} \in \mathbb{R}^{d^{\prime} \times d^{\prime}}$. Define the data-misfit function
\[
\Phi_{\bfy} (\bfSampPi) = \frac{1}{2}\left \Vert \bfSigma_{\bfepsilon}^{-\frac{1}{2}}\left(G(\bfSampPi) - \bfy\right)\right \Vert_2^2\,,
\]
with the Euclidean norm $\|\cdot \|_2$. The likelihood function $L_{\bfy}: \Pcal \to \mathbb{R}$ 
is 
\[
L_{\bfy}( \bfSampPi) = \exp\left(-\Phi_{\bfy}(\bfSampPi)\right)\,.
\]
Combing the prior distribution with density $\pi_0$ and the likelihood $L$ via Bayes' theorem gives the posterior density up to a normalizing constant
\[
\post(\bfSampPi) \propto L_{\bfy} ( \bfSampPi) \post_0(\bfSampPi)\,.
\]
Above and for the remainder of this paper, we assume that the prior and posterior measures are absolutely continuous with respect to the Lebesgue measure.

\subsection{The Metropolis-Hastings algorithm}
\label{sec:Prelim:MH}
MCMC methods are widely used to sample posterior distributions that arise in Bayesian inverse problems. The Metropolis-Hastings algorithm defines a wide class of MCMC methods, on which we will build in the following. Algorithm~\ref{alg:MCMC-MH} describes the Metropolis-Hastings approach. In our Bayesian inverse problem setting, inputs are the likelihood $L$, the prior density $\post_0$, a proposal density $q$, and the number of iterations $M \in \mathbb{N}$. In each iteration, a proposal sample $\bfSampPi^{\prime}$ is drawn from the proposal distribution, which may depend on the previous sample $\bfSampPi_{i - 1}$. Then, the acceptance probability $\alpha(\bfSampPi_{i - 1}, \bfSampPi^{\prime})$ is computed, which requires evaluating the likelihood $L$ and prior $\post_0$ at the candidate sample $\bfSampPi^{\prime}$. The proposal sample $\bfSampPi^{\prime}$ is accepted with probability $\alpha(\bfSampPi_{i - 1}, \bfSampPi^{\prime})$ and rejected with probability $1 - \alpha(\bfSampPi_{i - 1}, \bfSampPi^{\prime})$.  This process is repeated for $M$ iterations and the samples $\bfSampPi_1, \dots, \bfSampPi_M$ are returned. The corresponding Markov chain is, by construction, reversible for $\post$ and thus has $\post$ as a stationary distribution. With some relatively simple additional conditions on $\post$ and $q$, one can show that the chain converges to $\post$, from any starting point; see \cite{Roberts} for a full discussion. 

\begin{algorithm}[t]
\caption{Metropolis-Hastings}\label{alg:MCMC-MH}
\begin{algorithmic}[1]
\Procedure{MetropolisHastings}{$L, \post_0, q, M$}
\State Choose a starting point $\bfSampPi_0$
\For{$i = 1, \dots, M$}
\State Draw candidate $\bfSampPi^{\prime}$ from proposal $q(\cdot | \bfSampPi_{i - 1})$
\State Compute acceptance probability
\[
\alpha(\bfSampPi_{i - 1}, \bfSampPi^{\prime}) = \min\left\{1, \frac{q(\bfSampPi_{i - 1} | \bfSampPi^{\prime})L_{\bfy}(\bfSampPi^{\prime})\pi_0(\bfSampPi^{\prime})}{q(\bfSampPi^{\prime}|\bfSampPi_{i - 1})L_{\bfy}(\bfSampPi_{i - 1})\pi_0(\bfSampPi_{i-1})}\right\}
\]
\State Set the sample $\bfSampPi_i$ to
\[
\bfSampPi_i = \begin{cases}
\bfSampPi^{\prime}\,,\qquad &\text{with probability }\alpha(\bfSampPi_{i - 1}, \bfSampPi^{\prime})\,,\\
\bfSampPi_{i - 1}\,,\qquad &\text{with probability }1 - \alpha(\bfSampPi_{i - 1}, \bfSampPi^{\prime})
\end{cases}
\]
\EndFor
\State \Return $\bfSampPi_1, \dots, \bfSampPi_{M}$
\EndProcedure
\end{algorithmic}
\end{algorithm}

In general, the samples produced by MCMC are correlated; this correlation inflates the variance of any expectations estimated with MCMC samples, relative to an expectation estimated with uncorrelated Monte Carlo sample sets of the same size. The efficiency of an MCMC sampler can thus be measured with the \emph{effective sample size} (ESS) of any sample set it produces, which is inversely proportional to the \emph{integrated autocorrelation time} of the chain \cite{MCBookLiu}. To define the ESS, consider a function $f: \Pcal \to \mathbb{R}$ that is measurable with respect to the Lebesgue measure, and let us assume we are interested in estimating the expected value
\[
\mathbb{E}[f] = \int_{\Pcal} f(\bfSampPi)\pi(\bfSampPi) \, \mathrm d\bfSampPi\,,
\]
with respect to the posterior distribution $\pi$. Consider now the Monte Carlo estimator of $\mathbb{E}[f]$ that uses $n \in \mathbb{N}$ samples $ \{ \bfSampPi_i \}_{i=1}^{n}$,
\[
E_nf = \frac{1}{n}\sum_{i = 1}^n f(\bfSampPi_i)\,.
\]
The ESS of $ \{ f(\bfSampPi_i) \}_{i=1}^{n}$ is $n^* \in \mathbb{R}$ such that
\[
\Var[E_n f] = \frac{\Var[f]}{n^*}\,,
\]
where $\Var[E_n f]$ is the variance of the estimator $E_n f$ and $\Var[f]$ is the variance of $f(\bfSampPi)$ for $\bfSampPi \sim \pi$. In other words, $n^* \leq n$ is the number of independent Monte Carlo samples from $\pi$ that would be required to obtain an estimator with the same variance as $E_n f$. It can be shown that $n^* = n / \tau$, where $\tau$ is the integrated autocorrelation time associated with the chain $\{f(\bfSampPi_t)\}_t$. Better MCMC mixing corresponds to  smaller $\tau$ and larger ESS.

\subsection{Transport maps}
\label{sec:Prelim:TMap}
The use of transport maps in the context of Bayesian inference was introduced in \cite{ELMOSELHY20127815}. In particular, \cite{ELMOSELHY20127815} proposed a variational Bayesian approach involving transport. Rather than using importance sampling or MCMC to characterize the posterior distribution, this approach seeks a transport map that pushes forward a tractable ``reference'' distribution to the posterior, such that samples drawn from the reference and acted on by the map are distributed according to the posterior. Below we follow \cite{Marzouk2016} to introduce the notion of transport maps.

\subsubsection{Definition of transport maps}
Let $\mu_{\post}$ and $\mu_{\refP}$ be two probability measures on $\mathbb{R}^d$ that are absolutely continuous with respect to the Lebesgue measure. In the following, $\mu_{\post}$ is the target measure, which corresponds to the posterior distribution in our case, and $\mu_{\refP}$ is the reference measure, which typically is a Gaussian or another distribution from which we can draw independent samples efficiently. The probability density function corresponding to the target measure $\mu_{\post}$ is the posterior $\post$, and the probability density function corresponding to $\mu_{\refP}$ is denoted by $\refP$. A transport map is a function $\mapT: \mathbb{R}^d \to \mathbb{R}^d$ that pushes forward the reference $\mu_{\refP}$ to the target $\mu_{\post}$, which we write as
\begin{equation}
\mu_{\post} = \mapT_{\sharp}\mu_{\refP}\,,
\label{eq:Maps:TMapCondition}
\end{equation}
and which means that for any Borel set $A \subseteq \mathbb{R}^d$, it holds $\mu_{\post}(A) = \mu_{\refP}(T^{-1}(A))$. Existence of such maps is ensured by the absolute continuity of the reference and the target measure. There may be infinitely many transport maps that push forward a given reference to the target of interest. Uniqueness can be enforced by introducing a cost function that is minimized while imposing the constraint \eqref{eq:Maps:TMapCondition}. This construction leads to the notion of \emph{optimal transport}; see, e.g., \cite{Vershik2013,Villani2003,Villani2009}.

Instead of introducing a cost function to regularize the problem of finding a transport map, we directly impose structure on the map $T$ as in \cite{ELMOSELHY20127815,ParnoTMap,Marzouk2016}. In particular, we will seek lower triangular maps that are 
monotone increasing. The lower triangular structure of $T$ is as follows
\begin{equation}
T(\sampR_1, \sampR_2, \dots, \sampR_d) = \left[\begin{array}{l}
T_1(\sampR_1) \\
T_2(\sampR_1, \sampR_2)\\
T_3(\sampR_1, \sampR_2, \sampR_3)\\
\vdots\\
T_d(\sampR_1, \sampR_2, \dots, \sampR_d)
\end{array}\right]\,,
\label{eq:Prelim:TriangularStructure}
\end{equation}
where $\sampR_i$ denotes the $i$th component of $\bfSampR = [\sampR_1, \dots, \sampR_d]^T \in \mathbb{R}^{d}$ and where $T_i: \mathbb{R}^i \to \mathbb{R}$ is the $i$th component function of the map $T$. Monotonicity in this context corresponds to the condition that $T_i$ is a monotone increasing function of $\sampR_i$, for all $i=1 \ldots d$. This condition ensures that $\nabla T \succeq 0$ and $\det \nabla T \geq 0$; see \cite{Marzouk2016,2017arXiv170306131S} for more detail.
Since we assume that the reference and the target measures are absolutely continuous, existence and uniqueness of such a  lower-triangular transport map is guaranteed; this map is in fact the Knothe--Rosenblatt rearrangement \cite{doi:10.1137/120874850,doi:10.1137/080740647,rosenblatt1952}. 

\subsubsection{Numerical approximations of transport maps}
\label{sec:Prelim:MapsNumerics}
Following \cite{ELMOSELHY20127815}, we will obtain numerical approximations of the Knothe--Rosenblatt rearrangement by solving an optimization problem. Let $\Tcal$ be a finite-dimensional subspace of the space of all smooth lower triangular functions \eqref{eq:Prelim:TriangularStructure} from $\mathbb{R}^d$ into $\mathbb{R}^d$. Then, an approximation $\mapTApprox \in \Tcal$ of a transport map $T$ can be obtained via numerical optimization over the coefficients of the representation of $\mapTApprox$ in a basis of $\Tcal$. To set up the optimization problem, consider the \emph{pullback} $(\mapTApprox^{-1})_\sharp \mu_{\post} \equiv \mapTApprox^{\sharp}\mu_{\post}$ of $\mu_{\post}$ through a map $\mapTApprox$; the density of this pullback measure can be written as 
\begin{equation}
\pb(\bfSampR) = \pi(\mapTApprox(\bfSampR))|\det \nabla \mapTApprox(\bfSampR)|\,,
\label{eq:Prelim:Pullback}
\end{equation}
where $|\det \nabla \mapTApprox(\bfSampR)|$ is the absolute value of the determinant of the Jacobian $\nabla \mapTApprox(\bfSampR)$ of $\mapTApprox$ at $\bfSampR$. Note that the functions in $\Tcal$ are smooth in the sense that $\nabla \mapTApprox$ exists and is sufficiently regular; see, e.g., \cite{Marzouk2016,2017arXiv170306131S}.
Let  
\[
D_{\text{KL}}(\pi_1 || \pi_2) = \mathbb{E}_{\pi_1}\left(\log \frac{\pi_1}{\pi_2}\right)
\]
denote the Kullback--Leibler (KL) divergence of $\pi_1$ from $\pi_2$ (where $\pi_1$ and $\pi_2$ in the $\log$ term are densities). Then, a solution $\mapTApprox^* \in \Tcal$ of the optimization problem
\begin{equation}
\begin{aligned}
\min_{\mapTApprox \in \mathcal{T}} ~~~& D_{\text{KL}}(\refP||\pb)\,,\\
\text{s.t.}~ & \nabla \mapTApprox \succ 0\,,
\end{aligned}
\label{eq:Maps:OptiProblem}
\end{equation}
is an approximation of a transport map that pushes forward the reference measure $\mu_{\refP}$ to the target $\mu_{\post}$. The constraint $\nabla \mapTApprox \succ 0$ means that the Jacobian of $\mapTApprox$ is positive definite. If the approximation space $\Tcal$ is sufficiently rich such that $D_{\text{KL}}(\refP || \pb) = 0$, then we have $\mapTApprox^*{}^{\sharp}\mu_{\post} = \mu_{\refP}$ and $\mapTApprox^*_{\sharp}\mu_{\refP} = \mu_{\post}$\cite{Marzouk2016}. 

The KL divergence is not symmetric. The direction of the KL divergence here is chosen such that the expected value is taken with respect to the reference $\mu_{\refP}$, which is selected so that it can easily be sampled. Furthermore, the objective in \eqref{eq:Maps:OptiProblem} can be minimized without knowledge of the normalizing constant of $\post$. To see this, transform the objective $D_{\text{KL}}(\refP||\pb)$ into 
\begin{equation}
\begin{aligned}
D_{\text{KL}}(\refP || \pb) & = \mathbb{E}_{\refP}\left[\log\left(\frac{\refP}{\pb}\right)\right]\\
& = \mathbb{E}_{\refP} \left[\log \refP - \log \post \circ \mapTApprox - \log | \det \nabla \mapTApprox | \right]\,,
\end{aligned}
\label{eq:Prelim:ObjTransform}
\end{equation}
where we used the definition of $\pb$ in \eqref{eq:Prelim:Pullback}. The expectation $\mathbb{E}_{\refP}[\log \refP]$ is independent of the map $\mapTApprox$ and therefore can be ignored when minimizing the objective of \eqref{eq:Maps:OptiProblem}. Similarly, the normalizing constant of $\post$ leads to a constant term in \eqref{eq:Prelim:ObjTransform} that is independent of $\mapTApprox$ and therefore unnormalized evaluations of $\post$ are sufficient to optimize \eqref{eq:Maps:OptiProblem}. 

\subsection{Problem description}
\label{sec:Prelim:Problem}
We identify two challenges of directly relying on a solution $\mapTApprox^*$ of \eqref{eq:Maps:OptiProblem} to solve a Bayesian inverse problem, as proposed in \cite{ELMOSELHY20127815}. First, the KL divergence objective \eqref{eq:Prelim:ObjTransform}, which contains an expectation with respect the the reference measure $\mu_{\refP}$, is typically estimated with a Monte Carlo method because no closed form expression is available for general $\post$. Thus, the (unnormalized version of the) pullback density \eqref{eq:Prelim:Pullback} must be evaluated at a potentially large number of samples for each optimization iteration. Since each evaluation of the pullback entails an evaluation of the (unnormalized) posterior distribution $\post$, and thus of the forward model $G$, the optimization can become computationally intractable if $G$ is expensive to evaluate. Second, the push forward $\mapTApprox^*_{\sharp}\mu_{\refP}$ is only an approximation of the posterior $\mu_{\post}$. In particular, approximation errors may follow from the choice of the finite-dimensional approximation space $\mathcal{T}$, the finite number of Monte Carlo samples used to discretize the expectation in the objective, and any other errors in the numerical optimization. While this error can be estimated (see \cite{ELMOSELHY20127815,Marzouk2016}), reducing this error to an arbitrarily small threshold---e.g., by enriching $\mathcal{T}$---can be computationally expensive.

\section{Multifidelity preconditioned Metropolis-Hastings}
\label{sec:MFTMap}
We propose a multifidelity preconditioned Metropolis-Hastings (MFMH). 
Similar to the approach of Moselhy~\cite{ELMOSELHY20127815}, our multifidelity approach consists of an offline and an online phase. In the offline phase, a transport map is constructed rapidly by using a low-cost, low-fidelity approximation of the high-fidelity model. 
Then, in the online phase, a proposal distribution is derived from the transport map and used to sample the high-fidelity posterior via a modified Metropolis-Hastings algorithm, following the work of Parno \cite{ParnoTMap}. 
If the low-fidelity model is an accurate approximation of the high-fidelity model, and if the transport map captures the essential structure of $\pi$,
then the ESS of this MFMH will be higher and fewer online evaluations of the high-fidelity model will be required to achieve a given accuracy. 
In any case, any error associated with the transport map is corrected by Metropolization in the online phase; samples are thus drawn (asymptotically) from the posterior distribution corresponding to the high-fidelity model $G$. In other words, the MFMH approach offers the same convergence guarantees as standard MCMC.

\subsection{Approximation spaces and numerical optimization}

This section provides details on the numerical construction of transport maps from the reference measure to the target measure.

\subsubsection{Integrated squared parametrization}
Each component function $\mapTApprox_i, i = 1, \dots, d$ of the approximation $\mapTApprox$ of the map $T$ defined in \eqref{eq:Prelim:TriangularStructure} is parameterized with the integrated-squared ansatz \cite{Bigoni}, which is similar to the integrated-exponential ansatz introduced in \cite{2017arXiv170306131S}. Component function $\mapTApprox_i$ is parameterized as 
\begin{multline}
\mapTApprox_i(\sampR_1, \dots, \sampR_i; \bfbeta_i) \\= \mapTApprox_i^{(L)}(\sampR_1, \dots, \sampR_{i-1}; \bfbeta_i^{(L)}) + \int_{0}^{\sampR_i} \left( \mapTApprox_i^{(R)}(\sampR_1, \dots, \sampR_{i - 1}, t; \bfbeta^{(R)}_i)\right)^2 \mathrm dt\,,
\label{eq:MFTMap:IntegratedSquaredParam}
\end{multline}
where 
\begin{equation}
\bfbeta_i = \begin{bmatrix} \bfbeta^{(L)}_i & \bfbeta^{(R)}_i\end{bmatrix}
\label{eq:MFTMap:BetaParameter}
\end{equation}
is a parameter vector and $\mapTApprox_i^{(L)}: \mathbb{R}^{i - 1} \to \mathbb{R}$ and $\mapTApprox_i^{(R)}: \mathbb{R}^i \to \mathbb{R}$ are functions that are parameterized by the parameters $\bfbeta^{(L)}_i$ and $\bfbeta^{(R)}_i$, respectively.  The integrated-squared parameterization \eqref{eq:MFTMap:IntegratedSquaredParam} guarantees that the map $\mapTApprox$ is monotone and therefore automatically satisfies the constraint $\nabla \mapTApprox \succ 0$ in \eqref{eq:Maps:OptiProblem} \cite{2017arXiv170306131S}.

\subsubsection{Approximation space}
Following \cite{ParnoTMap,Marzouk2016}, we represent the functions $\mapTApprox_i^{(L)}$ and $\mapTApprox_i^{(R)}$ in each component $\mapTApprox_i, i = 1, \dots, d$ of the map $\mapTApprox$ as multivariate polynomials. Let $\bfj = [j_1, \dots, j_d]^T \in \mathbb{N}^d$ be a multi-index and let $\phi_{j_i}$ be a univariate polynomial with degree $j_i$ for $i = 1, \dots, d$. Define the multivariate polynomial function as
\begin{equation}
\phi_{\bfj}(\bfSampR) = \prod_{i = 1}^d \phi_{j_i}(\sampR_i)\,.
\label{eq:MFTMap:Mono}
\end{equation}
In the following, the multivariate polynomial functions \eqref{eq:MFTMap:Mono} are simply the monomials. Note that other polynomial families can be used, e.g., Hermite polynomials. (One might also use Hermite functions, as in \cite{2017arXiv170306131S}, for better control of tail behavior.) Consider now the sets $\Jcal_i \subset \mathbb{N}^d$ for $i = 1, \dots, d$
\begin{equation}
\Jcal_i = \{\bfj \,|\, \|\bfj\|_1 \leq \ell\,,\,\, j_k = 0\,,\,\, \forall k > i\}\,, 
\label{eq:MFTMap:Jcal}
\end{equation}
which correspond to the total-degree polynomials of maximal degree $\ell \in \mathbb{N}$. The constraint $j_k = 0\,,\, \forall k > i$ in the definition of $\Jcal_i$ in \eqref{eq:MFTMap:Jcal} imposes the lower-triangular structure of the map $\mapTApprox$ as defined in \eqref{eq:Prelim:TriangularStructure}. The set $\Jcal_i$ leads to the definition of the approximation space $\Tcal_i$
\[
\Tcal_i = \operatorname{span}\{\phi_{\bfj} \,|\, \bfj \in \Jcal_i\}\,,
\]
of the functions $\mapTApprox_i^{(R)}$ and $\mapTApprox_{i+1}^{(L)}$, respectively, of the $i$th and $i+1$st component of $\mapTApprox$. Thus, $\mapTApprox_i^{(L)} \in \Tcal_{i-1}$ and $\mapTApprox_i^{(R)} \in \Tcal_i$ can be represented as
\begin{equation}
\mapTApprox_i^{(L)}(\sampR_1, \dots, \sampR_{i-1}; \bfbeta_i^{(L)})  = \sum_{\bfj \in \Jcal_{i-1}} \beta_{i,\bfj}^{(L)} \phi_{\bfj}(\sampR_1, \dots, \sampR_{i-1})\,,
\label{eq:MFTMap:SparameterizationL}
\end{equation}
and
\begin{equation}
\mapTApprox_i^{(R)}(\sampR_1, \dots, \sampR_{i-1}, t; \bfbeta_i^{(R)}) = \sum_{\bfj \in \Jcal_i} \beta_{i,\bfj}^{(R)} \phi_{\bfj}(\sampR_1, \dots, \sampR_{i-1},t)\,,
\label{eq:MFTMap:SparameterizationR}
\end{equation}
where $\bfbeta_i^{(L)}$ and $\bfbeta_i^{(R)}$ are the vectors of the coefficients $\bfbeta_{i,\bfj}^{(L)}$ with $\bfj \in \Jcal_{i-1}$ and $\bfbeta_{i,\bfj}^{(R)}$ with $\bfj \in \Jcal_i$, respectively. Note that we combine $\bfbeta_i^{(L)}$ and $\bfbeta_i^{(R)}$ into the vector $\bfbeta_i$ as in defined in \eqref{eq:MFTMap:BetaParameter}. The approximation space $\Tcal$ is the product space
\[
\Tcal = (\Tcal_0 \oplus \Tcal_1) \otimes \dots \otimes (\Tcal_{d - 1} \oplus \Tcal_d)\,,
\]
where $\oplus$ and $\otimes$ denote the sum and the product of two spaces.

\subsubsection{Numerical solution of optimization problem}
\label{sec:MFMH:OptimizationProblem}
Consider the optimization problem \eqref{eq:Maps:OptiProblem} and the transformation \eqref{eq:Prelim:ObjTransform} that shows that it is sufficient to minimize 
\begin{equation}
\mathbb{E}_{\refP} \left[ - \log \pi \circ \mapTApprox - \log \det \nabla \mapTApprox \right]
\label{eq:Maps:ObjToOptimize}
\end{equation}
with respect to $\mapTApprox$. Note that evaluations of the unnormalized version of $\pi$ are sufficient to minimize \eqref{eq:Maps:ObjToOptimize}; see Section~\ref{sec:Prelim:MapsNumerics}. We have dropped the absolute value above since the maps $\mapTApprox$ are guaranteed to be monotone, via our parameterization. We now replace the expected value with its sample-average approximation \cite{SpallStochOpti}, i.e., a Monte Carlo estimator employing independent draws $ \{ \bfSampR_i \}_{i=1}^n$ from the reference distribution $\refP$. Making the dependence on the coefficients $\bfbeta_1, \dots, \bfbeta_d$ explicit, we obtain the optimization problem
\begin{equation}
\min_{\bfbeta_1, \dots, \bfbeta_d}  \frac{1}{n} \sum_{i = 1}^n \left[- \log \pi(\mapTApprox(\bfSampR_i; \bfbeta_1, \dots, \bfbeta_d)) - \log \det \nabla \mapTApprox(\bfSampR_i; \bfbeta_1, \dots, \bfbeta_d) \right]\,.
\label{eq:MFTMap:MFOptiProblem}
\end{equation}
Our optimization problem \eqref{eq:MFTMap:MFOptiProblem} is unconstrained because the constraint $\nabla \mapTApprox \succ 0$ of \eqref{eq:Maps:OptiProblem} is automatically satisfied via the squared-integrated parameterization \eqref{eq:MFTMap:IntegratedSquaredParam}.

\subsection{Constructing transport maps with low-fidelity models}
We propose to approximate the high-fidelity forward model $G$ with a low-cost, low-fidelity model to reduce the computational costs of constructing a transport map via the optimization problem \eqref{eq:MFTMap:MFOptiProblem}. Note that replacing evaluations of $G$ with evaluations of low-fidelity models in \eqref{eq:MFTMap:MFOptiProblem} will introduce an error that we must correct later; see Section~\ref{sec:MFTMap:Precond}.

Let $\widehat{G}: \Pcal \to \Ycal$ be a low-fidelity approximation of $G$ that maps the parameter onto an observable. The model $\widehat{G}$ gives rise to a low-fidelity potential function
\[
\widehat{\Phi}_{\bfy}(\bfSampPi) = \frac{1}{2}\left\|\bfSigma_{\bfepsilon}^{-\frac{1}{2}}\left(\widehat{G}(\bfSampPi) - \bfy\right)\right\|_2^2\,,
\]
and to the low-fidelity likelihood function $\widehat{L}_{\bfy}: \Pcal \to \mathbb{R}$
\[
\widehat{L}_{\bfy}(\bfSampPi) = \exp\left(-\widehat{\Phi}_{\bfy}(\bfSampPi)\right)\,.
\]
The probability density of the corresponding low-fidelity posterior distribution is, up to a normalizing constant,
\[
\postLow(\bfSampPi) \propto \widehat{L}_{\bfy} (\bfSampPi)\post_0(\bfSampPi)\,.
\]

We now use the low-fidelity posterior density $\postLow$ in our variational construction of the map $\mapTApproxLow \in \Tcal$. 
Consider the density $\pbLow$ that is the pullback of the low-fidelity posterior density $\postLow$ through a map $\mapTApproxLow$, 
\[
\pbLow(\bfSampR) = \postLow(\mapTApproxLow(\bfSampR))|\det \nabla \mapTApproxLow(\bfSampR)|\,.
\]
We can find $\mapTApproxLow$ by minimizing $D_{\text{KL}}(\refP || \pbLow)$; this yields the following optimization problem, analogous to \eqref{eq:MFTMap:MFOptiProblem} but with $\postLow$ replacing $\post$,
\begin{equation}
\min_{\bfbeta_1, \dots, \bfbeta_d}  \frac{1}{n} \sum_{i = 1}^n \left[- \log \postLow(\mapTApproxLow(\bfSampR_i; \bfbeta_1, \dots, \bfbeta_d)) - \log \det \nabla \mapTApproxLow(\bfSampR_i; \bfbeta_1, \dots, \bfbeta_d) \right]\,,
\label{eq:MFTMap:MFOptiProblemLowFid}
\end{equation}
with the coefficients $\bfbeta_1, \dots, \bfbeta_d$ defining $\mapTApproxLow^* \in \Tcal$.

\subsection{Constructing deep transport maps}
\label{sec:MFMH:Compose}
The optimization problem \eqref{eq:MFTMap:MFOptiProblemLowFid} finds a transport map in the approximation space $\Tcal$. Thus, if the approximation space $\Tcal$ is chosen too coarse, then the pullback of the target distribution through the map is only a poor approximation of the reference distribution. Instead of choosing richer approximation spaces to find more accurate transport maps, which would lead to a large number of coefficients to be optimized for in \eqref{eq:MFTMap:MFOptiProblemLowFid}, the work \cite{ELMOSELHY20127815,ParnoThesis}
proposes to take compositions of transport maps in coarse approximation spaces.

Let $\mapT^{(1)}$ be a transport map that pushes forward the reference $\mu_{\refP}$ onto the posterior $\mu_{\postLow}$ and let $\mapTApproxLow^{(1)} \in \Tcal$ be a numerical approximation of $\mapT^{(1)}$ derived with optimization problem \eqref{eq:MFTMap:MFOptiProblemLowFid}. Since $\mapTApproxLow^{(1)}$ is an approximation of the map $\mapT^{(1)}$, 
we obtain 
\[
\mu_{\postLow} \approx \mu_{\postLow}^{(1)} = \mapTApproxLow^{(1)}_{\sharp} \mu_{\refP}\,.
\]
Following \cite{ParnoTMap}, to account for the discrepancy between $\mu_{\postLow}$ and $\mu_{\postLow}^{(1)}$, a second map $\mapTApproxLow^{(2)} \in \Tcal$ is constructed to push forward $\mu_{\postLow}^{(1)}$ to $\mu_{\postLow}$, i.e., 
\[
\mu_{\postLow} \approx \mu_{\postLow}^{(2)} = \mapTApproxLow^{(2)}_{\sharp} \mu_{\postLow}^{(1)}\,.
\]
The aim is that the map $\mapTApproxLow^{(2)}$ should capture only a small correction from $\mu^{(1)}_{\postLow}$ to $\mu^{(2)}_{\postLow}$, so that the composition $\mapTApproxLow^{(1, 2)} = \mapTApproxLow^{(2)} \circ \mapTApproxLow^{(1)}$ more accurately pushes forward $\mu_{\refP}$ to $\mu_{\postLow}$ (in the sense of KL divergence). This process is  repeated $k$ times to obtain the ``deep'' map
\begin{equation}
\mapTApproxLow^* = \mapTApproxLow^{(k)} \circ \dots \circ \mapTApproxLow^{(1)}\,.
\label{eq:TMap:Composition}
\end{equation}

Note that \eqref{eq:TMap:Composition} is only one possible way of constructing deep transport maps. Another possibility is to reverse the order of \eqref{eq:TMap:Composition} and compose maps from the right, rather than from the left, by keeping the $\mu_{\refP}$ as the reference and pulling back the target $\mu_{\postLow}$ by the current map approximation. This construction is discussed in \cite{2017arXiv170306131S}.

\subsection{MFMH algorithm}
\label{sec:MFTMap:Precond}
We now use the transport map $\mapTApproxLow^*$ to derive an MCMC proposal with the aim of improving sampling efficiency (e.g., reducing integrated autocorrelation time) and at the same time guaranteeing that the stationary distribution of the Markov chain is the posterior distribution corresponding to the high-fidelity model $G$. 

Since the transport map $\mapTApproxLow^*$ is monotone and lower triangular, as enforced by the approximation space $\Tcal$, the inverse map $\mapTInvApproxLow^* = \mapTApproxLow^*{}^{-1}$ is cheap to evaluate by solving $d$ one-dimensional root-finding problems. We refer to \cite{ParnoTMap,Marzouk2016} for details on how to invert lower-triangular, monotone maps. The inverse map $\mapTInvApproxLow^*$ pushes forward the true posterior $\mu_{\post}$ onto a distribution that \emph{approximates} the reference distribution. In other words, $\mapTInvApproxLow^*_\sharp \mu_{\post} = \mu_{\refP}$ only if the map were exact, which in general it is not; otherwise, $\mapTInvApproxLow^*_\sharp \mu_{\post}$ is simply closer to $\mu_{\refP}$ (in the sense of KL divergence) than $\mu_{\post}$ was. If $\mu_{\refP}$ is chosen to be a Gaussian, then the map $\mapTInvApproxLow^*$ approximately ``Gaussianizes'' $\mu_{\post}$. As shown in \cite{ParnoTMap}, and as will be demonstrated with our numerical results in Section~\ref{sec:NumExp}, MCMC sampling from the approximate reference distribution $\mapTInvApproxLow^*_\sharp \mu_{\post}$ generally yields higher ESSs for a given computational effort than sampling directly from the posterior distribution. 
Crucially, because the map $\mapTInvApproxLow^* = \mapTApproxLow^*{}^{-1}$ is invertible, samples from the approximate reference distribution can be \emph{exactly} (up to machine precision) pushed forward to samples from the posterior distribution via the map $\mapTApproxLow^*$.

Our MFMH algorithm is summarized in Algorithm~\ref{alg:PreMCMC-MH} and follows the preconditioned Metropolis-Hastings algorithm introduced in \cite{ParnoTMap}. The MFMH algorithm has the same inputs as the single-fidelity Metropolis-Hastings in Algorithm~\ref{alg:MCMC-MH}, except that additionally the maps $\mapTApproxLow^*$ and $\mapTInvApproxLow^*$ are required. The current state $\bfSampPi_{i - 1}$ is mapped with $\mapTInvApproxLow^*$ onto the approximate reference to obtain $\bfSampR_{i - 1}$. Then, a candidate sample is drawn from the proposal distribution based on $\bfSampR_{i - 1}$, and that candidate sample is mapped back with $\mapTApproxLow^*$ to a sample of the posterior $\pi$. 

Whether the candidate sample is accepted or rejected is based on the probability $\alpha$ that is determined using the high-fidelity posterior density $\pi$; the maps, derived from the low-fidelity posterior, only affect the proposal distribution. Thus it is guaranteed that the stationary distribution of the chain is the posterior distribution $\pi$. Note that there is considerable flexibility in the choice of the reference-space proposal distribution $q$. Below we will mostly choose $q(\cdot \vert \bfSampR_{i-1} )$ to be a $d$-dimensional Gaussian, \emph{independent} of $\bfSampR_{i-1}$. Thus our MCMC algorithm reduces to a Metropolis independence sampler \cite{Roberts}, whether viewed on the reference space or on the target space. Given the monotonicity of the transport maps and the full support of the Gaussian reference on $\mathbb{R}^d$, the proposal distribution on the target space is guaranteed to dominate the posterior distribution, as required for Metropolis independence sampling to converge.

\begin{algorithm}[t]
\caption{Multifidelity preconditioned Metropolis-Hastings (MFMH)}\label{alg:PreMCMC-MH}
\begin{algorithmic}[1]
\Procedure{PreconditionedMetropolisHastings}{$L_{\bfy}, \pi_0, q, M, \mapTApproxLow^*, \mapTInvApproxLow^*$}
\State Choose a starting point $\bfSampPi_0$
\For{$i = 1, \dots, M$}
\State Map state $\bfSampPi_{i - 1}$ onto reference with $\mapTInvApproxLow^*(\bfSampPi_{i - 1}) = \bfSampR_{i - 1}$
\State Draw candidate $\bfSampR^{\prime}$ from proposal $q(\cdot | \bfSampR_{i - 1})$
\State Map candidate $\bfSampR^{\prime}$ onto target with $\mapTApproxLow^*(\bfSampR^{\prime}) = \bfSampPi^{\prime}$
\State Compute acceptance probability
\[
\alpha(\bfSampPi_{i - 1}, \bfSampPi^{\prime}) = \min\left\{1, \frac{q(\bfSampR_{i - 1} |\bfSampR^{\prime})L_{\bfy}(\bfSampPi^{\prime})\pi_0(\bfSampPi^{\prime})|\det \nabla \mapTApproxLow^*(\bfSampPi^{\prime})|}{q(\bfSampR^{\prime}|\bfSampR_{i - 1})L_{\bfy}(\bfSampPi_{i - 1})\pi_0(\bfSampPi_{i-1})|\det \nabla \mapTApproxLow^*(\bfSampPi_{i - 1})|}\right\}
\]
\State Set the sample $\bfSampPi_i$ to
\[
\bfSampPi_i = \begin{cases}
\bfSampPi^{\prime}\,,\qquad &\text{with probability }\alpha(\bfSampPi_{i - 1}, \bfSampPi^{\prime})\,,\\
\bfSampPi_{i - 1}\,,\qquad &\text{with probability }1 - \alpha(\bfSampPi_{i - 1}, \bfSampPi^{\prime})
\end{cases}
\]
\EndFor
\State \Return $\bfSampPi_1, \dots, \bfSampPi_{M}$
\EndProcedure
\end{algorithmic}
\end{algorithm}

\section{Numerical results}
\label{sec:NumExp}
This section demonstrates our multifidelity approach on two examples. All runtime measurements were performed on compute nodes with Intel Xeon E5-1620 and 64GB RAM on a single core using a {\scshape Matlab} implementation.

\subsection{Diffusion equation with reaction term}
We first consider a model with diffusion and a nonlinear reaction term, where our goal is to infer the parameters of the reaction term. 

\subsubsection{Problem setup}
Let $\Omega = (0, 1)^2 \subseteq \mathbb{R}^2$ and $\Pcal = \mathbb{R}^2$ and consider the PDE
\begin{equation}
-\nabla^2 u(x_1, x_2; \bfSampPi) + g(u(x_1, x_2; \bfSampPi), \bfSampPi) = 100\sin(2\pi x_1)\sin(2\pi x_2)\,, \qquad \bfx \in \Omega\,,
\label{eq:NumExp:Laplace:PDE}
\end{equation}
with homogeneous Dirichlet boundary conditions, where $\bfx = [x_1, x_2]^T$, $\bfSampPi = [\sampPi_1, \sampPi_2]^T \in \Pcal  = \mathbb{R}^2$, and $u: \Omega \times \Pcal \to \mathbb{R}$ is the solution function. The nonlinear function $g$ is
\[
g(u(\bfx; \bfSampPi), \bfSampPi) = (0.1\sin(\sampPi_1) + 2)\exp\left(-2.7\sampPi_1^2 \right)\left(\exp\left(1.8\sampPi_2 u(\bfx; \bfSampPi) \right) - 1\right)\,.
\] 
We discretize \eqref{eq:NumExp:Laplace:PDE} with finite differences on a grid with equidistant grid points and mesh width $h > 0$. The corresponding system of nonlinear equations is solved with Newton's method and inexact line search based on the Armijo condition. The model $G_h: \Pcal \to \Ycal$ derived with mesh width $h$ maps from $\Pcal$ into $\Ycal = \mathbb{R}^{12}$. The components of the observable $\bfy$ correspond to the values of the approximated solution function at the spatial coordinates $[0.25i, 0.2j]^T \in \Omega$ with $i = 1, 2, 3$ and $j = 1, 2, 3, 4$. 

A low-fidelity model $\widehat{G}_h$ of $G_h$ is derived via projection-based model reduction \cite{BennerSIREV}. Solutions of \eqref{eq:NumExp:Laplace:PDE} for parameters on an $100 \times 100$ equidistant grid in $[-\pi/2, -\pi/2] \times [1, 5] \subset \Pcal$ and mesh width $h$ are computed and a 20-dimensional reduced space with proper orthogonal decomposition is constructed. The operators corresponding to the high-fidelity $G_h$ are projected via Galerkin projection onto the reduced space and the low-fidelity model $\widehat{G}_h$ is obtained. 

\subsubsection{Setup of inverse problems}
\label{sec:NumExp:Laplace:InvProb}
We set $\bfSampPi^* = [0.5, 2]^T$ and consider the data $\bfy = G_H(\bfSampPi^*) + \bfepsilon$, where $H = 1/64$ and $\bfepsilon$ adds Gaussian noise with zero mean and variance $0.0026$, which corresponds to $0.1 \%$ noise with respect to the Euclidean norm of the output $G_H(\bfSampPi^*)$. We then use the model $G_h$ with $h = 1/32$ and the corresponding low-fidelity model $\widehat{G}_h$ for inference. Note that the data $\bfy$ is computed with a discretization of the PDE \eqref{eq:NumExp:Laplace:PDE} with mesh width $H = 1/64$, while we use a mesh width $h = 1/32$ for inference. The prior distribution is a Gaussian distribution with mean $[\pi/4, 1.2]^T$ and covariance matrix
\[
\bfSigma_{\bfepsilon} = \begin{bmatrix}
1 & 0\\
0 & 0.01
\end{bmatrix}\,.
\]
The low-fidelity $\widehat{G}_h$ is about 80 times faster to evaluate than the high-fidelity model $\widehat{G}_h$ for $h = 1/32$.

We construct a transport map from the reference Gaussian distribution to the posterior corresponding to the low-fidelity model $\widehat{G}_h$. The reference Gaussian distribution has zero mean and standard deviation $0.1$. Note that instead we could have scaled and centered the posterior distribution to have zero mean and marginal variances of one and then used a standard Gaussian as reference distribution.

The optimization is performed with {\scshape Matlab}'s \texttt{fmincon} optimizer, where the tolerance (\texttt{TolX}) is set to $10^{-3}$ and where we use $n = 250$ samples of the reference distribution to approximate the expected value in the objective function, see Section~\ref{sec:MFMH:OptimizationProblem}. We compose two maps following Section~\ref{sec:MFMH:Compose}. The approximation space for the first map $\mapTApproxLow^{(1)}$ corresponds to first-order polynomials (linear) and the approximation space of the second map $\mapTApproxLow^{(2)}$ corresponds to second-order polynomials (quadratic). The starting point for the optimization is the identity map. The transport map is then used to precondition Metropolis-Hastings as shown in Algorithm~\ref{alg:PreMCMC-MH}. The proposal is the reference distribution, i.e., we obtain an independence sampler with a proposal that is independent of the current state of the chain. Because we expect that the transport map approximately pushes the reference onto the posterior, it is reasonable to consider an independence sampler with the reference distribution as proposal. We discard every other sample, which means that we perform $M = 2m$ iterations if we want $m$ samples. 

We compare sampling with our MFMH approach to delayed-rejection adaptive-metropolis (DRAM) sampling \cite{Haario2006}. We initialize DRAM with a Gaussian proposal that has a diagonal covariance with all elements on the diagonal being equal. We start DRAM with diagonal elements in $\{10^{-4}, 5 \times 10^{-4}, 10^{-3}, 5 \times 10^{-3}, 10^{-2},2 \times 10^{-2}, 3 \times 10^{-2}, 4 \times 10^{-2}, 5 \times 10^{-2}, 10^{-1}, 5 \times 10^{-1}\}$ and then select the run with the highest ESS, cf.~Section~\ref{sec:Prelim:MH}. The first $10^4$ samples are discarded as burn-in and then every other sample is used. This means that DRAM performs $M = 2m + 10^4$ iterations if we want $m$ samples. Note that the same thinning of discarding every other sample is applied to the samples obtained with MFMH and DRAM. 

\begin{figure}
\begin{tabular}{cc}
{\huge\resizebox{0.48\columnwidth}{!}{\input{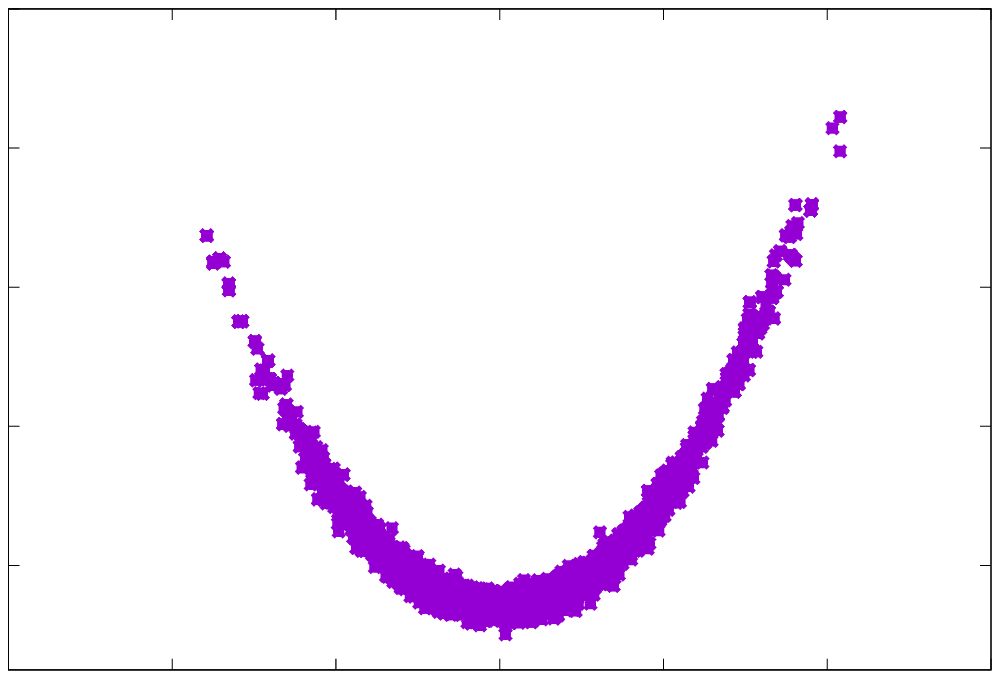}}} & {{\huge\resizebox{0.48\columnwidth}{!}{\input{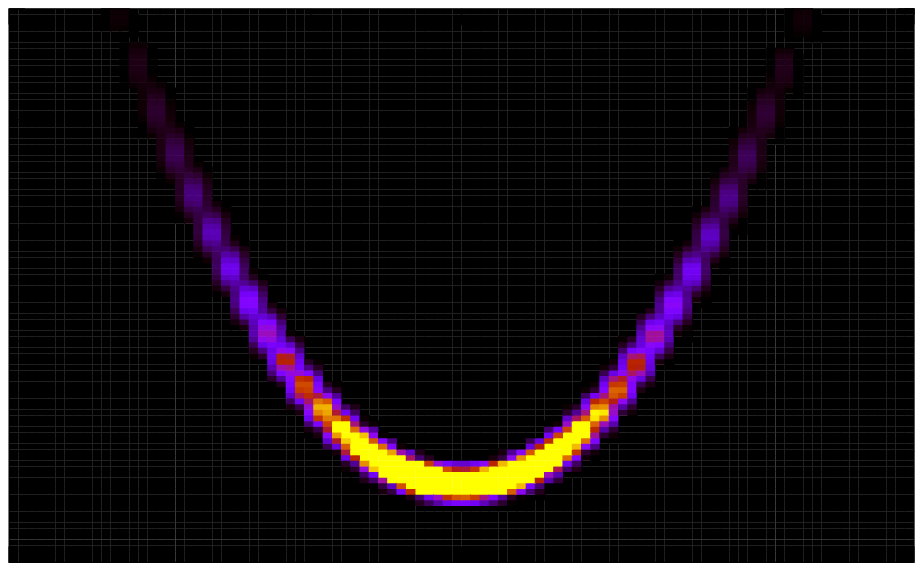}}}}\\
\scriptsize (a) samples from DRAM & \scriptsize (b) posterior corresponding to high-fidelity model\\
{\huge\resizebox{0.48\columnwidth}{!}{\input{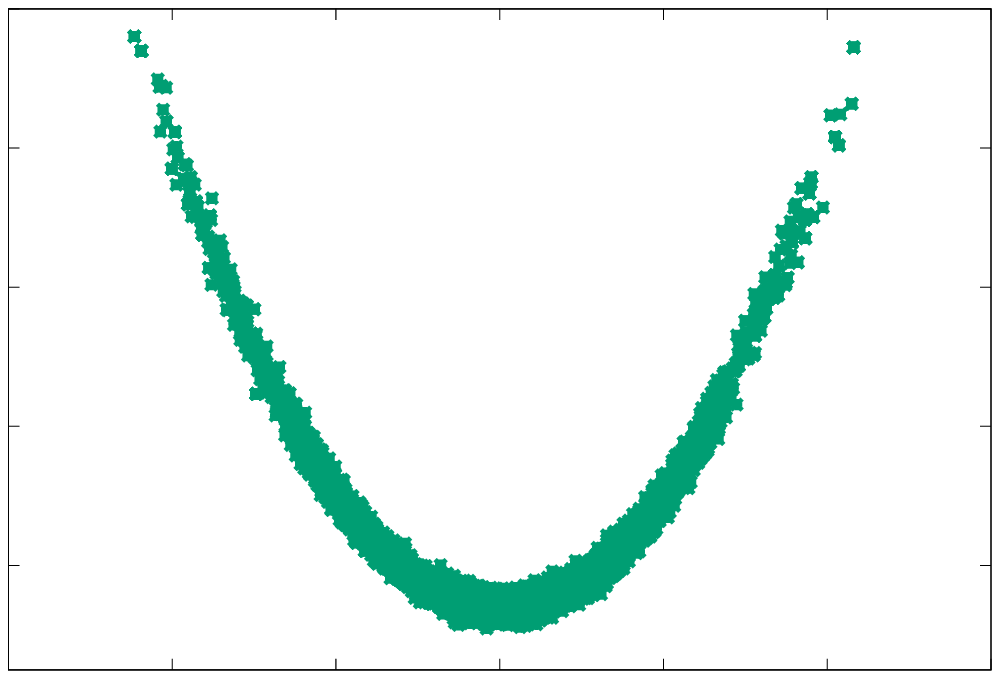}}} & {{\huge\resizebox{0.48\columnwidth}{!}{\input{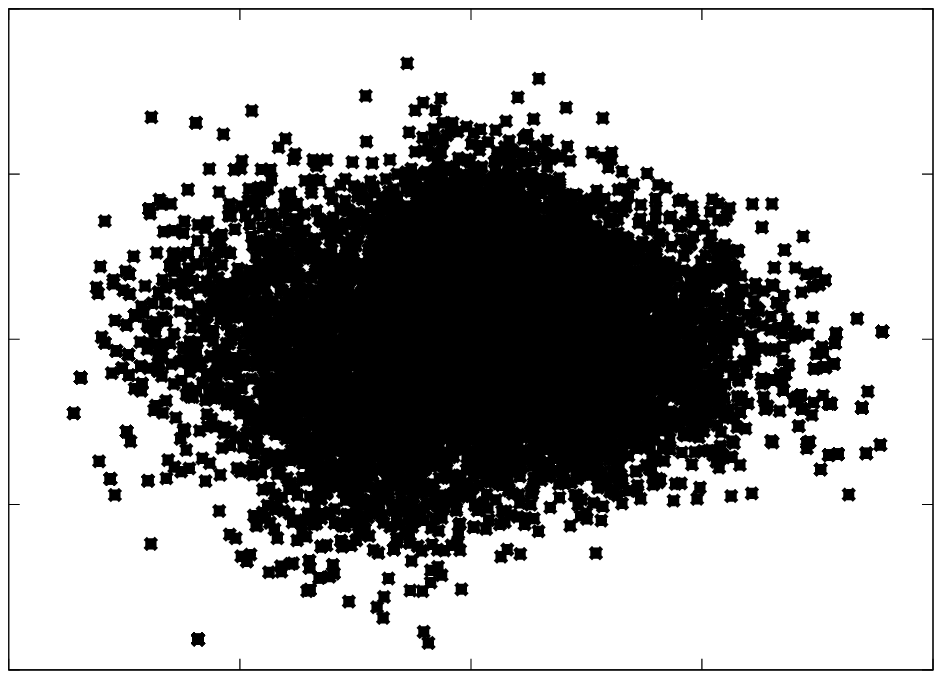}}}}\\
 \scriptsize (c) samples from MFMH & {\scriptsize (d) samples from MFMH mapped onto reference}
\end{tabular}
\caption{Diffusion-reaction problem: Plots in (a) and (c) visualize samples drawn with DRAM and MFMH, respectively. The ``banana''-like shape of the posterior is reflected in both sets of samples, see plot (b). Plot (d) shows the samples from (c) mapped onto the reference distribution, which demonstrates that our transport map captures well the high-fidelity posterior distribution even though it was constructed with a low-fidelity model.}
\label{fig:Laplace:Data}
\end{figure}

\subsubsection{Results}
\begin{figure}
\begin{center}
{\large\resizebox{0.7\columnwidth}{!}{\input{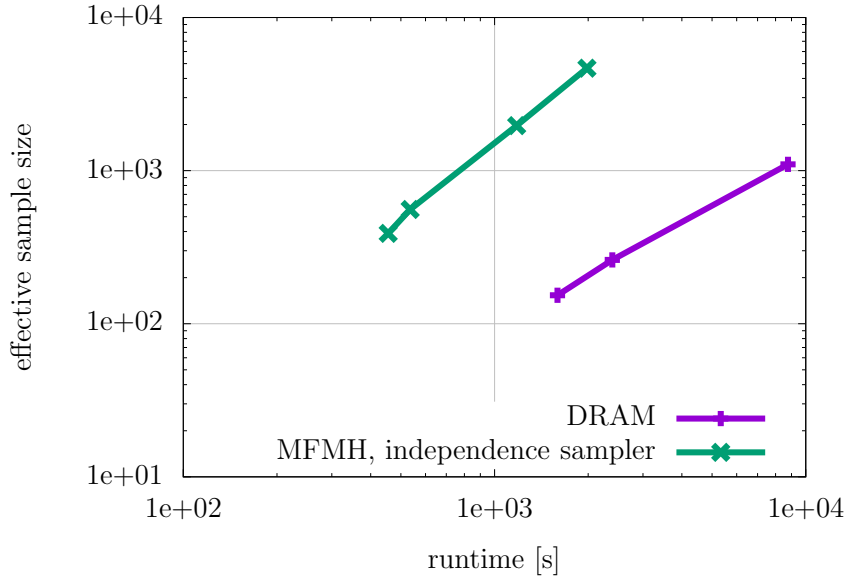}}}
\end{center}
\caption{Diffusion-reaction problem: The ESS of samples drawn with our MFMH approach is higher than with DRAM, in this example. The reported runtime includes the runtime of constructing the transport map in case of our MFMH approach. Note that there is no burn-in time for MFMH because it is based on an independence sampler in this example.}
\label{fig:Laplace:ESS}
\end{figure}

Figure~\ref{fig:Laplace:Data} shows samples drawn with DRAM from the high-fidelity posterior and compares them to samples drawn with MFMH. Both reflect the ``banana''-like shape of the posterior. Mapping the samples from Figure~\ref{fig:Laplace:Data}c to the reference distribution results in the samples shown in Figure~\ref{fig:Laplace:Data}d. Figure~\ref{fig:Laplace:ESS} reports the ESS of samples drawn with DRAM and our MFMH approach. For MFMH, we report results for $m \in \{5 \times 10^2, 10^3, 5 \times 10^3, 10^4\}$ samples and for DRAM for $m \in \{5 \times 10^3, 10^4, 5 \times 10^4\}$ samples. In this example, MFMH achieves a higher ESS than DRAM with respect to runtime. The reported runtime is the time needed for model evaluations and, in case of MFMH, for constructing and evaluating the transport map. Note that we use a burn-in of $10^4$ samples for DRAM, whereas such a burn-in is unnecessary in case of our MFMH sampler because it is an independence sampler.

\subsection{Euler Bernoulli beam problem}
\label{sec:NumExp:Euler}
We now infer the effective stiffness of an Euler Bernoulli beam, for which a model is available on GitHub\footnote{https://github.com/g2s3-2018/labs}. The model was developed by Matthew Parno for the 2018 Gene Golub SIAM Summer School on ``Inverse Problems: Systematic Integration of Data with Models under Uncertainty.'' 

\subsubsection{Problem setup}
\label{sec:NumExp:Euler:Setup}
Let $L > 0$ be the length of the beam and define $\Omega = [0, L] \subset \mathbb{R}$. Consider the PDE
\begin{equation}
\frac{\partial^2}{\partial x^2}\left(E(x)\frac{\partial^2}{\partial x^2}u(x)\right) = f(x)\,,\qquad x \in \bar{\Omega}\,,
\label{eq:Euler:PDE}
\end{equation}
where $u: \Omega \to \mathbb{R}$ is the vertical deflection of the beam, $f: \Omega \to \mathbb{R}$ is the load, and $\bar{\Omega} = (0, L)$. The effective stiffness of the beam is given by $E: \Omega \to \mathbb{R}$ and describes beam geometry and material properties. The beam is in cantilever configuration, where the left boundary is fixed and the right boundary is free, i.e., the boundary conditions are
\[
u(0) = 0\,,\quad \left.\frac{\partial}{\partial x}u\right|_{x = 0} = 0\,,\quad \left.\frac{\partial^2}{\partial x^2} u\right|_{x = L} = 0\,,\quad \left.\frac{\partial^3}{\partial x^3}u\right|_{x = L} = 0\,.
\]
The length of the beam is $L=1$ in the following. The PDE \eqref{eq:Euler:PDE} is discretized with finite differences on a mesh of $N = 601$ equidistant grid points in $\Omega$. The same effective stiffness $E$ as available in the GitHub repository\footnotemark[1] is used by interpolating on the grid points. Let $x_{\text{obs}}^{(1)}, \dots, x_{\text{obs}}^{(41)}$ be equidistant points of the $N = 601$ grid points. The observation $\bfy \in \mathbb{R}^{41}$ is the displacement $u$ at the 41 points $x_{\text{obs}}^{(1)}, \dots, x_{\text{obs}}^{(41)}$ polluted with zero-mean Gaussian noise with variance $10^{-4}$.

\subsubsection{Models for the Euler Bernoulli problem}
We now derive the high-fidelity forward model. Consider the function $I: \mathbb{R} \times \Omega \to \mathbb{R}$ defined as
\[
I(x, \alpha) = \left(1 + \exp\left(-\frac{x - \alpha}{0.005}\right)\right)^{-1}\,,
\]
with 
\[
\lim_{x \to -\infty} I(x, \alpha) = 0\,,\qquad \lim_{x \to \infty}I(x, \alpha) = 1\,,
\]
such that there is a smooth transition from $0$ to $1$ at $\alpha$. Define further $k = 3$ and let $\alpha_1, \dots, \alpha_{k+1}$ be the $k + 1$ equidistant points in $\Omega$. Let $\mathbb{R}_+ = \{z \in \mathbb{R} \,:\, z > 0\}$ and consider the parameter $\bfSampPi = [\sampPi_1, \dots, \sampPi_k]^T \in \mathbb{R}^k_+$. Define the function $\widehat{E}_i: \Omega \times \mathbb{R} \to \mathbb{R}$ as
\[
\widehat{E}_i(x, \sampPi_i) = (1 - I(x, \alpha_i))\widehat{E}_{i-1}(x, \sampPi_{i-1}) + I(x, \alpha_i)\sampPi_i\,,
\]
for $i = 2, \dots, k$ and $\widehat{E}_1(x, \sampPi_1) = \sampPi_1$. Given a parameter $\bfSampPi$, the function $\widehat{E}_k$ is a smooth approximation of the piecewise constant function $\sum_{i = 1}^k \sampPi_i \mathbb{I}_{(\alpha_i, \alpha_{i+1}]}$, where $\mathbb{I}_{(\alpha_i, \alpha_{i + 1}]}$ is the indicator function of the interval $(\alpha_i, \alpha_{i + 1}] \subset \mathbb{R}$. The high-fidelity forward model $G$ maps a parameter $\bfSampPi \in \mathbb{R}_{+}^k$ onto the displacement $u$ with effective stiffness $\widehat{E}_k$ at the observation points $x_{\text{obs}}^{(1)}, \dots, x_{\text{obs}}^{(41)}$. The same discretization as described in Section~\ref{sec:NumExp:Euler:Setup} is used. The low-fidelity model $\widehat{G}$ is a spline interpolant of $G$ on a logarithmically spaced grid in the domain $[0.5, 4]^3$. The interpolant is obtained with the \texttt{griddedInterpolant} procedure available in {\scshape Matlab}. Extrapolation is turned off. The low-fidelity model $\widehat{G}$ is about $1,100$ times faster to evaluate than the high-fidelity model $G$. 

\begin{figure}
\begin{tabular}{c}
{\large\resizebox{0.7\columnwidth}{!}{\input{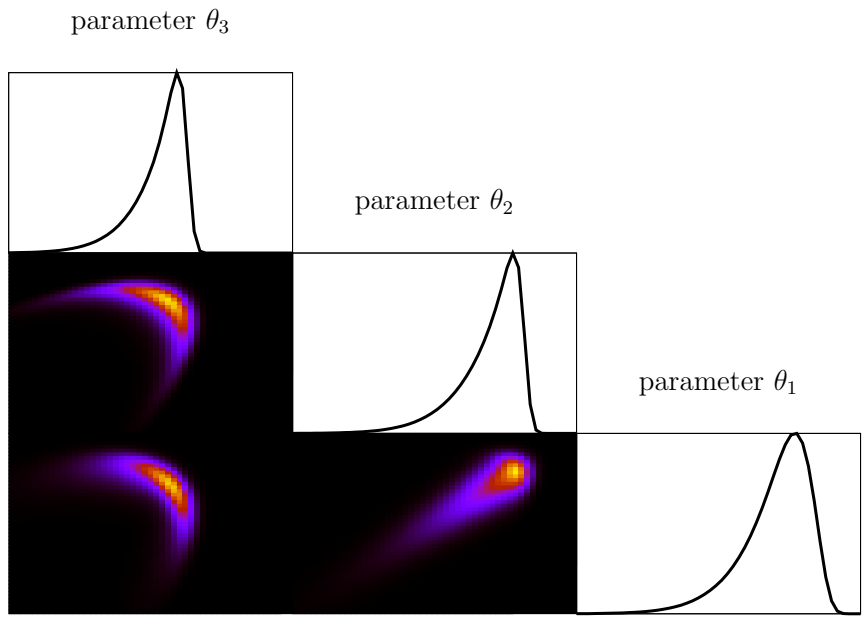}}}\\
(a) dimension-wise heat map of posterior density function\\
~\\
~\\
{\large\resizebox{0.7\columnwidth}{!}{\input{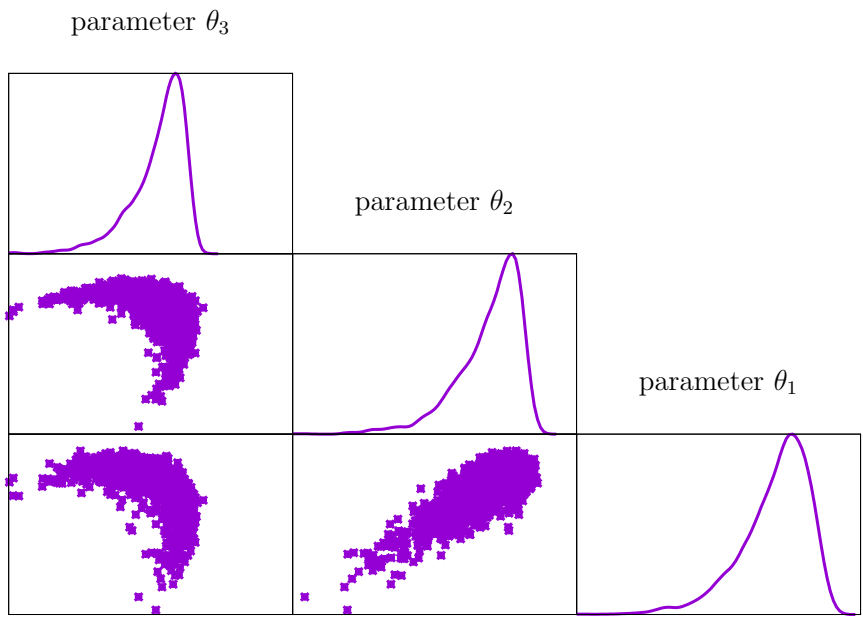}}}\\
(b) samples drawn from posterior distribution with DRAM 
\end{tabular}
\caption{Euler Bernoulli beam: Plot (a) visualizes the posterior density function. Plot (b) shows samples drawn with DRAM from the posterior distribution. The plots on the diagonal of (b) show marginal densities estimated from the samples.}
\label{fig:Euler:DRAM}
\end{figure}

\subsubsection{Setup of inverse problem}
The observation $\bfy$ is obtained as described in Section~\ref{sec:NumExp:Euler:Setup}. The prior is a log-normal distribution with mean $[1, 1, 1]^T$ and covariance matrix
\[
\bfSigma_{\bfepsilon} =  \begin{bmatrix}
0.05 & 0 & 0\\
0 & 0.05 & 0\\
0 & 0 & 0.05
\end{bmatrix}\,.
\]
We construct a transport map $\mapTApproxLow$ using the low-fidelity model $\widehat{G}$. The transport map $\mapTApproxLow$ is based on quadratic polynomials and approximately maps the reference Gaussian distribution with mean $[1, 1, 1]^T$ and covariance with diagonal entries $0.1$ onto the low-fidelity posterior distribution. Note that a standard Gaussian can be used as reference distribution if the posterior distribution is centered to zero mean and scaled to have marginal variances of one. This effectively would perform a linear transformation. 
We optimize for the coefficients of $\mapTApproxLow$ with \texttt{fmincon} available in {\scshape Matlab} with tolerance set to $10^{-8}$ (\texttt{TolX} option) and with 500 samples from the reference distribution. The initial point of the optimization corresponds to the identity map. We then use the transport map in Algorithm~\ref{alg:PreMCMC-MH} with a random-walk Metropolis algorithm on the reference space, i.e., with a local Gaussian proposal $\bfSampR^{\prime} \sim \mathcal{N}(\bfSampR_{i - 1}, \bfSigma^{\prime})$. The covariance matrix $\bfSigma^{\prime}$ is diagonal. To select the diagonal of $\bfSigma^{\prime}$, we ran MFMH for proposals with variance in $\{10^{-4}, 5 \times 10^{-4}, 10^{-3}, 5 \times 10^{-3}, 10^{-2},2 \times 10^{-2}, 3 \times 10^{-2}, 4 \times 10^{-2}, 5 \times 10^{-2}, 10^{-1}, 5 \times 10^{-1}\}$ and then selected the variance that leads to the highest ESS, cf.~Section~\ref{sec:NumExp:Laplace:InvProb}. The corresponding results are reported in the following.
Additionally, we consider an independence MFMH sampler where the reference distribution together with the transport map serves as a proposal distribution that is independent of the previous sample. The rest of the setup is the same as in Section~\ref{sec:NumExp:Laplace:InvProb}. 

\subsubsection{Results}
Figure~\ref{fig:Euler:DRAM}a visualizes the posterior density function corresponding to the high-fidelity model. The plot indicates that the posterior distribution is a non-Gaussian distribution. We use the same procedure as in Section~\ref{sec:NumExp:Laplace:InvProb} to draw samples from the posterior with DRAM. The samples with the highest ESS are shown in Figure~\ref{fig:Euler:DRAM}b. 

\begin{figure}
\begin{tabular}{c}
{\large\resizebox{0.7\columnwidth}{!}{\input{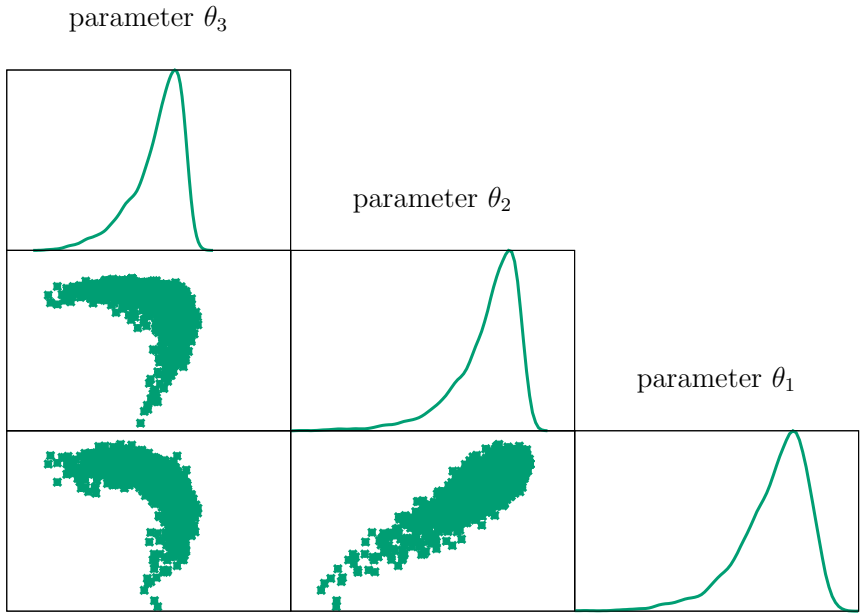}}}\\
(a) samples drawn with MFMH\\
~\\
~\\
{\large\resizebox{0.7\columnwidth}{!}{\input{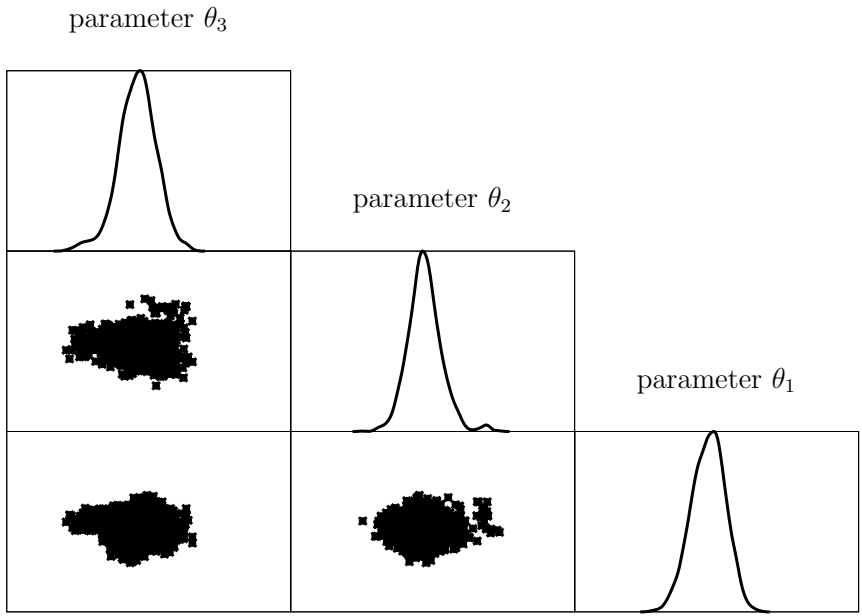}}}\\
(b) samples mapped onto an approximation of the reference Gaussian distribution
\end{tabular}
\caption{Euler Bernoulli beam: Plot (a) shows samples drawn with our MFMH method using a transport map constructed with a low-fidelity model. Mapping the samples shown in (a) with the transport map gives samples of an approximate reference Gaussian distribution, see plot (b). The plots on the diagonal of (a) and (b) show marginal densities estimated from the samples.}
\label{fig:Euler:MFMH}
\end{figure}

\begin{figure}
\begin{center}
{\resizebox{0.6\columnwidth}{!}{\input{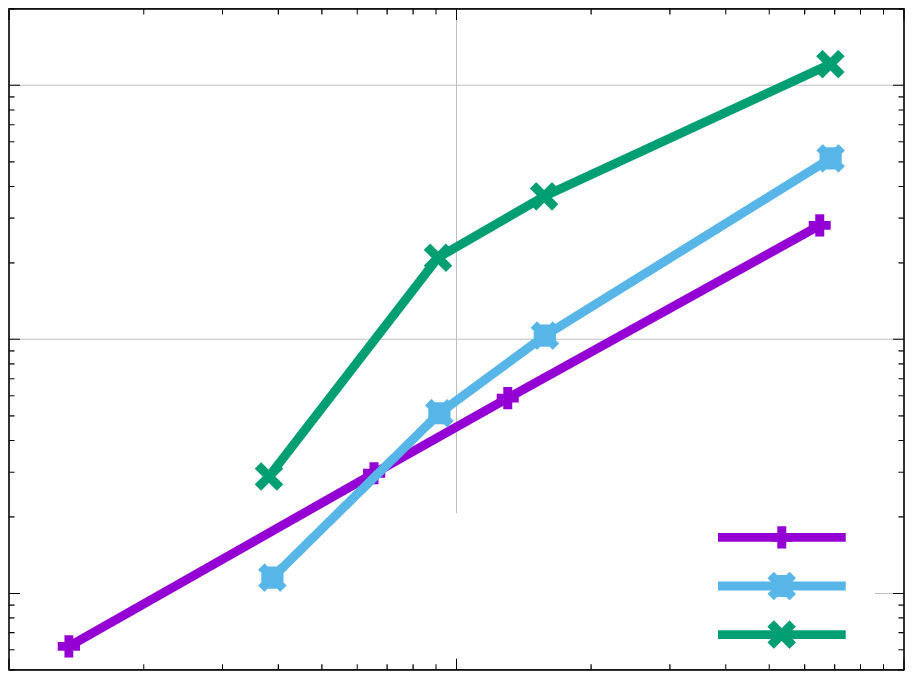}}}
\end{center}
\caption{Euler Bernoulli beam: Our MFMH approach leads to a higher ESS than DRAM after the offline costs of constructing the transport map are compensated. The runtime is the total time of drawing samples, including the runtime of constructing the transport map with the low-fidelity model in case of sampling with MFMH. The points on the curves correspond to $m \in \{10^5, 5 \times 10^5, 10^{6}, 5 \times 10^6\}$ samples, respectively.}
\label{fig:Euler:ESS}
\end{figure}

Let us now consider our MFMH approach. First, we consider the independence MFMH sampler that uses a proposal that is independent of the previous sample, see Figure~\ref{fig:Euler:MFMH}a. The samples are in agreement with the samples drawn with DRAM from the high-fidelity posterior distribution, cf.~Figure~\ref{fig:Euler:DRAM}b. Figure~\ref{fig:Euler:MFMH}b plots the samples mapped with the transport map onto an approximation of the reference Gaussian distribution, which demonstrates that the transport map $\mapTApproxLow^*$ captures well the posterior distribution. The ESSs of $m \in \{10^5, 5 \times 10^5, 10^{6}, 5 \times 10^6\}$ samples drawn with DRAM and the independence MFMH sampler are compared in Figure~\ref{fig:Euler:ESS}. Our independence MFMH sampler achieves a higher ESS than DRAM in this example. Figure~\ref{fig:Euler:ESS} also shows that a higher ESS compared to DRAM is achieved only if sufficiently many samples are drawn because otherwise the offline costs of constructing the transport maps are not compensated. The MFMH sampler with a local random walk proposal 
leads to a higher ESS than DRAM too, but the improvement is smaller compared to the independence MFMH sampler. Additionally, an even larger number of samples is necessary to compensate the offline costs of constructing the transport maps. Note that the MFMH sampler with a local random walk proposal and the DRAM sampler both use a burn-in of $10^4$ samples as in Section~\ref{sec:NumExp:Laplace:InvProb}. 

\section{Conclusions}
\label{sec:Conc}
Our MFMH approach leverages low-fidelity models to precondition MCMC sampling, with the aim of improving MCMC mixing while guaranteeing that the stationary distribution of the chain is the posterior distribution corresponding to the high-fidelity model. In a one-time expensive step, a low-fidelity model is used to construct a transport map that approximately maps an easy-to-sample reference distribution to the posterior distribution corresponding to the low-fidelity model. In the second step, the transport map is used to precondition the posterior distribution corresponding to the high-fidelity model during sampling with Metropolis-Hastings. Since the low-fidelity model is used for preconditioning only, the stationary distribution of the chain obtained in the second step is the posterior distribution corresponding to the high-fidelity model. Our MFMH approach achieves significant speedups compared to single-fidelity sampling with DRAM in our numerical examples.

\bibliography{multitm}
\bibliographystyle{abbrv}

\end{document}

%% file: figures/ExpDataDRAM.tex
\begingroup
  \makeatletter
  \providecommand\color[2][]{%
    \GenericError{(gnuplot) \space\space\space\@spaces}{%
      Package color not loaded in conjunction with
      terminal option `colourtext'%
    }{See the gnuplot documentation for explanation.%
    }{Either use 'blacktext' in gnuplot or load the package
      color.sty in LaTeX.}%
    \renewcommand\color[2][]{}%
  }%
  \providecommand\includegraphics[2][]{%
    \GenericError{(gnuplot) \space\space\space\@spaces}{%
      Package graphicx or graphics not loaded%
    }{See the gnuplot documentation for explanation.%
    }{The gnuplot epslatex terminal needs graphicx.sty or graphics.sty.}%
    \renewcommand\includegraphics[2][]{}%
  }%
  \providecommand\rotatebox[2]{#2}%
  \@ifundefined{ifGPcolor}{%
    \newif\ifGPcolor
    \GPcolortrue
  }{}%
  \@ifundefined{ifGPblacktext}{%
    \newif\ifGPblacktext
    \GPblacktexttrue
  }{}%
  \let\gplgaddtomacro\g@addto@macro
  \gdef\gplbacktext{}%
  \gdef\gplfronttext{}%
  \makeatother
  \ifGPblacktext
    \def\colorrgb#1{}%
    \def\colorgray#1{}%
  \else
    \ifGPcolor
      \def\colorrgb#1{\color[rgb]{#1}}%
      \def\colorgray#1{\color[gray]{#1}}%
      \expandafter\def\csname LTw\endcsname{\color{white}}%
      \expandafter\def\csname LTb\endcsname{\color{black}}%
      \expandafter\def\csname LTa\endcsname{\color{black}}%
      \expandafter\def\csname LT0\endcsname{\color[rgb]{1,0,0}}%
      \expandafter\def\csname LT1\endcsname{\color[rgb]{0,1,0}}%
      \expandafter\def\csname LT2\endcsname{\color[rgb]{0,0,1}}%
      \expandafter\def\csname LT3\endcsname{\color[rgb]{1,0,1}}%
      \expandafter\def\csname LT4\endcsname{\color[rgb]{0,1,1}}%
      \expandafter\def\csname LT5\endcsname{\color[rgb]{1,1,0}}%
      \expandafter\def\csname LT6\endcsname{\color[rgb]{0,0,0}}%
      \expandafter\def\csname LT7\endcsname{\color[rgb]{1,0.3,0}}%
      \expandafter\def\csname LT8\endcsname{\color[rgb]{0.5,0.5,0.5}}%
    \else
      \def\colorrgb#1{\color{black}}%
      \def\colorgray#1{\color[gray]{#1}}%
      \expandafter\def\csname LTw\endcsname{\color{white}}%
      \expandafter\def\csname LTb\endcsname{\color{black}}%
      \expandafter\def\csname LTa\endcsname{\color{black}}%
      \expandafter\def\csname LT0\endcsname{\color{black}}%
      \expandafter\def\csname LT1\endcsname{\color{black}}%
      \expandafter\def\csname LT2\endcsname{\color{black}}%
      \expandafter\def\csname LT3\endcsname{\color{black}}%
      \expandafter\def\csname LT4\endcsname{\color{black}}%
      \expandafter\def\csname LT5\endcsname{\color{black}}%
      \expandafter\def\csname LT6\endcsname{\color{black}}%
      \expandafter\def\csname LT7\endcsname{\color{black}}%
      \expandafter\def\csname LT8\endcsname{\color{black}}%
    \fi
  \fi
    \setlength{\unitlength}{0.0500bp}%
    \ifx\gptboxheight\undefined%
      \newlength{\gptboxheight}%
      \newlength{\gptboxwidth}%
      \newsavebox{\gptboxtext}%
    \fi%
    \setlength{\fboxrule}{0.5pt}%
    \setlength{\fboxsep}{1pt}%
\begin{picture}(7200.00,5040.00)%
    \gplgaddtomacro\gplbacktext{%
      \csname LTb\endcsname%
      \put(868,1497){\makebox(0,0)[r]{\strut{}$2$}}%
      \put(868,2299){\makebox(0,0)[r]{\strut{}$4$}}%
      \put(868,3100){\makebox(0,0)[r]{\strut{}$6$}}%
      \put(868,3902){\makebox(0,0)[r]{\strut{}$8$}}%
      \put(868,4703){\makebox(0,0)[r]{\strut{}$10$}}%
      \put(1036,616){\makebox(0,0){\strut{}$-3$}}%
      \put(1979,616){\makebox(0,0){\strut{}$-2$}}%
      \put(2922,616){\makebox(0,0){\strut{}$-1$}}%
      \put(3866,616){\makebox(0,0){\strut{}$0$}}%
      \put(4809,616){\makebox(0,0){\strut{}$1$}}%
      \put(5752,616){\makebox(0,0){\strut{}$2$}}%
      \put(6695,616){\makebox(0,0){\strut{}$3$}}%
    }%
    \gplgaddtomacro\gplfronttext{%
      \csname LTb\endcsname%
      \put(224,2799){\rotatebox{-270}{\makebox(0,0){\strut{}parameter $\theta_2$}}}%
      \put(3865,196){\makebox(0,0){\strut{}parameter $\theta_1$}}%
    }%
    \gplbacktext
    \put(0,0){\includegraphics{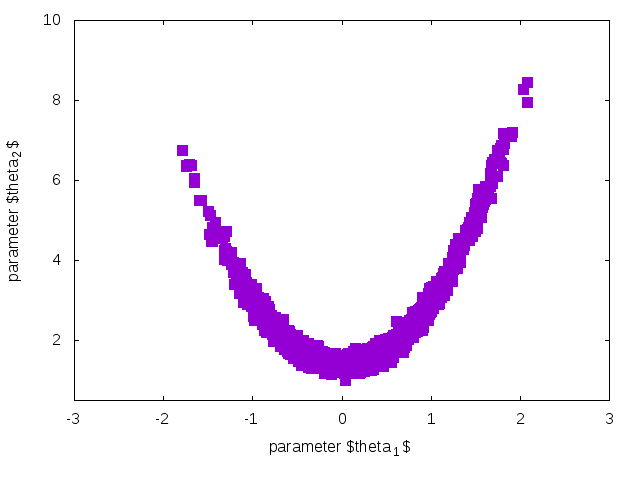}}%
    \gplfronttext
  \end{picture}%
\endgroup

%% file: figures/ExpPost.tex
\begingroup
  \makeatletter
  \providecommand\color[2][]{%
    \GenericError{(gnuplot) \space\space\space\@spaces}{%
      Package color not loaded in conjunction with
      terminal option `colourtext'%
    }{See the gnuplot documentation for explanation.%
    }{Either use 'blacktext' in gnuplot or load the package
      color.sty in LaTeX.}%
    \renewcommand\color[2][]{}%
  }%
  \providecommand\includegraphics[2][]{%
    \GenericError{(gnuplot) \space\space\space\@spaces}{%
      Package graphicx or graphics not loaded%
    }{See the gnuplot documentation for explanation.%
    }{The gnuplot epslatex terminal needs graphicx.sty or graphics.sty.}%
    \renewcommand\includegraphics[2][]{}%
  }%
  \providecommand\rotatebox[2]{#2}%
  \@ifundefined{ifGPcolor}{%
    \newif\ifGPcolor
    \GPcolortrue
  }{}%
  \@ifundefined{ifGPblacktext}{%
    \newif\ifGPblacktext
    \GPblacktexttrue
  }{}%
  \let\gplgaddtomacro\g@addto@macro
  \gdef\gplbacktext{}%
  \gdef\gplfronttext{}%
  \makeatother
  \ifGPblacktext
    \def\colorrgb#1{}%
    \def\colorgray#1{}%
  \else
    \ifGPcolor
      \def\colorrgb#1{\color[rgb]{#1}}%
      \def\colorgray#1{\color[gray]{#1}}%
      \expandafter\def\csname LTw\endcsname{\color{white}}%
      \expandafter\def\csname LTb\endcsname{\color{black}}%
      \expandafter\def\csname LTa\endcsname{\color{black}}%
      \expandafter\def\csname LT0\endcsname{\color[rgb]{1,0,0}}%
      \expandafter\def\csname LT1\endcsname{\color[rgb]{0,1,0}}%
      \expandafter\def\csname LT2\endcsname{\color[rgb]{0,0,1}}%
      \expandafter\def\csname LT3\endcsname{\color[rgb]{1,0,1}}%
      \expandafter\def\csname LT4\endcsname{\color[rgb]{0,1,1}}%
      \expandafter\def\csname LT5\endcsname{\color[rgb]{1,1,0}}%
      \expandafter\def\csname LT6\endcsname{\color[rgb]{0,0,0}}%
      \expandafter\def\csname LT7\endcsname{\color[rgb]{1,0.3,0}}%
      \expandafter\def\csname LT8\endcsname{\color[rgb]{0.5,0.5,0.5}}%
    \else
      \def\colorrgb#1{\color{black}}%
      \def\colorgray#1{\color[gray]{#1}}%
      \expandafter\def\csname LTw\endcsname{\color{white}}%
      \expandafter\def\csname LTb\endcsname{\color{black}}%
      \expandafter\def\csname LTa\endcsname{\color{black}}%
      \expandafter\def\csname LT0\endcsname{\color{black}}%
      \expandafter\def\csname LT1\endcsname{\color{black}}%
      \expandafter\def\csname LT2\endcsname{\color{black}}%
      \expandafter\def\csname LT3\endcsname{\color{black}}%
      \expandafter\def\csname LT4\endcsname{\color{black}}%
      \expandafter\def\csname LT5\endcsname{\color{black}}%
      \expandafter\def\csname LT6\endcsname{\color{black}}%
      \expandafter\def\csname LT7\endcsname{\color{black}}%
      \expandafter\def\csname LT8\endcsname{\color{black}}%
    \fi
  \fi
    \setlength{\unitlength}{0.0500bp}%
    \ifx\gptboxheight\undefined%
      \newlength{\gptboxheight}%
      \newlength{\gptboxwidth}%
      \newsavebox{\gptboxtext}%
    \fi%
    \setlength{\fboxrule}{0.5pt}%
    \setlength{\fboxsep}{1pt}%
\begin{picture}(7200.00,5040.00)%
    \gplgaddtomacro\gplbacktext{%
    }%
    \gplgaddtomacro\gplfronttext{%
      \csname LTb\endcsname%
      \put(994,667){\makebox(0,0){\strut{}$-3$}}%
      \put(1863,667){\makebox(0,0){\strut{}$-2$}}%
      \put(2732,667){\makebox(0,0){\strut{}$-1$}}%
      \put(3600,667){\makebox(0,0){\strut{}$0$}}%
      \put(4468,667){\makebox(0,0){\strut{}$1$}}%
      \put(5337,667){\makebox(0,0){\strut{}$2$}}%
      \put(6206,667){\makebox(0,0){\strut{}$3$}}%
      \put(3600,247){\makebox(0,0){\strut{}parameter $\theta_1$}}%
      \put(755,1704){\makebox(0,0)[r]{\strut{}$2$}}%
      \put(755,2342){\makebox(0,0)[r]{\strut{}$4$}}%
      \put(755,2978){\makebox(0,0)[r]{\strut{}$6$}}%
      \put(755,3616){\makebox(0,0)[r]{\strut{}$8$}}%
      \put(755,4254){\makebox(0,0)[r]{\strut{}$10$}}%
      \put(335,2660){\rotatebox{-270}{\makebox(0,0){\strut{}parameter $\theta_2$}}}%
    }%
    \gplbacktext
    \put(0,0){\includegraphics{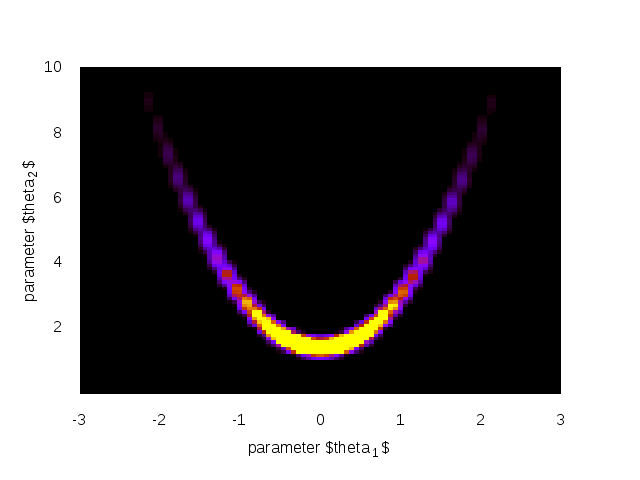}}%
    \gplfronttext
  \end{picture}%
\endgroup

%% file: figures/ExpDataPrecond.tex
\begingroup
  \makeatletter
  \providecommand\color[2][]{%
    \GenericError{(gnuplot) \space\space\space\@spaces}{%
      Package color not loaded in conjunction with
      terminal option `colourtext'%
    }{See the gnuplot documentation for explanation.%
    }{Either use 'blacktext' in gnuplot or load the package
      color.sty in LaTeX.}%
    \renewcommand\color[2][]{}%
  }%
  \providecommand\includegraphics[2][]{%
    \GenericError{(gnuplot) \space\space\space\@spaces}{%
      Package graphicx or graphics not loaded%
    }{See the gnuplot documentation for explanation.%
    }{The gnuplot epslatex terminal needs graphicx.sty or graphics.sty.}%
    \renewcommand\includegraphics[2][]{}%
  }%
  \providecommand\rotatebox[2]{#2}%
  \@ifundefined{ifGPcolor}{%
    \newif\ifGPcolor
    \GPcolortrue
  }{}%
  \@ifundefined{ifGPblacktext}{%
    \newif\ifGPblacktext
    \GPblacktexttrue
  }{}%
  \let\gplgaddtomacro\g@addto@macro
  \gdef\gplbacktext{}%
  \gdef\gplfronttext{}%
  \makeatother
  \ifGPblacktext
    \def\colorrgb#1{}%
    \def\colorgray#1{}%
  \else
    \ifGPcolor
      \def\colorrgb#1{\color[rgb]{#1}}%
      \def\colorgray#1{\color[gray]{#1}}%
      \expandafter\def\csname LTw\endcsname{\color{white}}%
      \expandafter\def\csname LTb\endcsname{\color{black}}%
      \expandafter\def\csname LTa\endcsname{\color{black}}%
      \expandafter\def\csname LT0\endcsname{\color[rgb]{1,0,0}}%
      \expandafter\def\csname LT1\endcsname{\color[rgb]{0,1,0}}%
      \expandafter\def\csname LT2\endcsname{\color[rgb]{0,0,1}}%
      \expandafter\def\csname LT3\endcsname{\color[rgb]{1,0,1}}%
      \expandafter\def\csname LT4\endcsname{\color[rgb]{0,1,1}}%
      \expandafter\def\csname LT5\endcsname{\color[rgb]{1,1,0}}%
      \expandafter\def\csname LT6\endcsname{\color[rgb]{0,0,0}}%
      \expandafter\def\csname LT7\endcsname{\color[rgb]{1,0.3,0}}%
      \expandafter\def\csname LT8\endcsname{\color[rgb]{0.5,0.5,0.5}}%
    \else
      \def\colorrgb#1{\color{black}}%
      \def\colorgray#1{\color[gray]{#1}}%
      \expandafter\def\csname LTw\endcsname{\color{white}}%
      \expandafter\def\csname LTb\endcsname{\color{black}}%
      \expandafter\def\csname LTa\endcsname{\color{black}}%
      \expandafter\def\csname LT0\endcsname{\color{black}}%
      \expandafter\def\csname LT1\endcsname{\color{black}}%
      \expandafter\def\csname LT2\endcsname{\color{black}}%
      \expandafter\def\csname LT3\endcsname{\color{black}}%
      \expandafter\def\csname LT4\endcsname{\color{black}}%
      \expandafter\def\csname LT5\endcsname{\color{black}}%
      \expandafter\def\csname LT6\endcsname{\color{black}}%
      \expandafter\def\csname LT7\endcsname{\color{black}}%
      \expandafter\def\csname LT8\endcsname{\color{black}}%
    \fi
  \fi
    \setlength{\unitlength}{0.0500bp}%
    \ifx\gptboxheight\undefined%
      \newlength{\gptboxheight}%
      \newlength{\gptboxwidth}%
      \newsavebox{\gptboxtext}%
    \fi%
    \setlength{\fboxrule}{0.5pt}%
    \setlength{\fboxsep}{1pt}%
\begin{picture}(7200.00,5040.00)%
    \gplgaddtomacro\gplbacktext{%
      \csname LTb\endcsname%
      \put(868,1497){\makebox(0,0)[r]{\strut{}$2$}}%
      \put(868,2299){\makebox(0,0)[r]{\strut{}$4$}}%
      \put(868,3100){\makebox(0,0)[r]{\strut{}$6$}}%
      \put(868,3902){\makebox(0,0)[r]{\strut{}$8$}}%
      \put(868,4703){\makebox(0,0)[r]{\strut{}$10$}}%
      \put(1036,616){\makebox(0,0){\strut{}$-3$}}%
      \put(1979,616){\makebox(0,0){\strut{}$-2$}}%
      \put(2922,616){\makebox(0,0){\strut{}$-1$}}%
      \put(3866,616){\makebox(0,0){\strut{}$0$}}%
      \put(4809,616){\makebox(0,0){\strut{}$1$}}%
      \put(5752,616){\makebox(0,0){\strut{}$2$}}%
      \put(6695,616){\makebox(0,0){\strut{}$3$}}%
    }%
    \gplgaddtomacro\gplfronttext{%
      \csname LTb\endcsname%
      \put(224,2799){\rotatebox{-270}{\makebox(0,0){\strut{}parameter $\theta_2$}}}%
      \put(3865,196){\makebox(0,0){\strut{}parameter $\theta_1$}}%
    }%
    \gplbacktext
    \put(0,0){\includegraphics{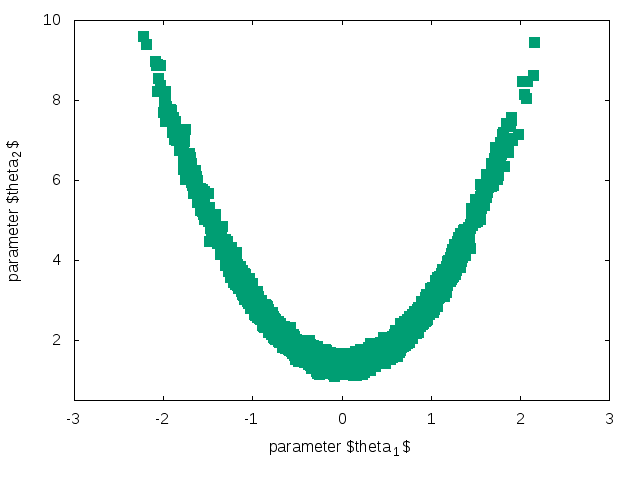}}%
    \gplfronttext
  \end{picture}%
\endgroup

%% file: figures/ExpDataPrecondProj.tex
\begingroup
  \makeatletter
  \providecommand\color[2][]{%
    \GenericError{(gnuplot) \space\space\space\@spaces}{%
      Package color not loaded in conjunction with
      terminal option `colourtext'%
    }{See the gnuplot documentation for explanation.%
    }{Either use 'blacktext' in gnuplot or load the package
      color.sty in LaTeX.}%
    \renewcommand\color[2][]{}%
  }%
  \providecommand\includegraphics[2][]{%
    \GenericError{(gnuplot) \space\space\space\@spaces}{%
      Package graphicx or graphics not loaded%
    }{See the gnuplot documentation for explanation.%
    }{The gnuplot epslatex terminal needs graphicx.sty or graphics.sty.}%
    \renewcommand\includegraphics[2][]{}%
  }%
  \providecommand\rotatebox[2]{#2}%
  \@ifundefined{ifGPcolor}{%
    \newif\ifGPcolor
    \GPcolortrue
  }{}%
  \@ifundefined{ifGPblacktext}{%
    \newif\ifGPblacktext
    \GPblacktexttrue
  }{}%
  \let\gplgaddtomacro\g@addto@macro
  \gdef\gplbacktext{}%
  \gdef\gplfronttext{}%
  \makeatother
  \ifGPblacktext
    \def\colorrgb#1{}%
    \def\colorgray#1{}%
  \else
    \ifGPcolor
      \def\colorrgb#1{\color[rgb]{#1}}%
      \def\colorgray#1{\color[gray]{#1}}%
      \expandafter\def\csname LTw\endcsname{\color{white}}%
      \expandafter\def\csname LTb\endcsname{\color{black}}%
      \expandafter\def\csname LTa\endcsname{\color{black}}%
      \expandafter\def\csname LT0\endcsname{\color[rgb]{1,0,0}}%
      \expandafter\def\csname LT1\endcsname{\color[rgb]{0,1,0}}%
      \expandafter\def\csname LT2\endcsname{\color[rgb]{0,0,1}}%
      \expandafter\def\csname LT3\endcsname{\color[rgb]{1,0,1}}%
      \expandafter\def\csname LT4\endcsname{\color[rgb]{0,1,1}}%
      \expandafter\def\csname LT5\endcsname{\color[rgb]{1,1,0}}%
      \expandafter\def\csname LT6\endcsname{\color[rgb]{0,0,0}}%
      \expandafter\def\csname LT7\endcsname{\color[rgb]{1,0.3,0}}%
      \expandafter\def\csname LT8\endcsname{\color[rgb]{0.5,0.5,0.5}}%
    \else
      \def\colorrgb#1{\color{black}}%
      \def\colorgray#1{\color[gray]{#1}}%
      \expandafter\def\csname LTw\endcsname{\color{white}}%
      \expandafter\def\csname LTb\endcsname{\color{black}}%
      \expandafter\def\csname LTa\endcsname{\color{black}}%
      \expandafter\def\csname LT0\endcsname{\color{black}}%
      \expandafter\def\csname LT1\endcsname{\color{black}}%
      \expandafter\def\csname LT2\endcsname{\color{black}}%
      \expandafter\def\csname LT3\endcsname{\color{black}}%
      \expandafter\def\csname LT4\endcsname{\color{black}}%
      \expandafter\def\csname LT5\endcsname{\color{black}}%
      \expandafter\def\csname LT6\endcsname{\color{black}}%
      \expandafter\def\csname LT7\endcsname{\color{black}}%
      \expandafter\def\csname LT8\endcsname{\color{black}}%
    \fi
  \fi
    \setlength{\unitlength}{0.0500bp}%
    \ifx\gptboxheight\undefined%
      \newlength{\gptboxheight}%
      \newlength{\gptboxwidth}%
      \newsavebox{\gptboxtext}%
    \fi%
    \setlength{\fboxrule}{0.5pt}%
    \setlength{\fboxsep}{1pt}%
\begin{picture}(7200.00,5040.00)%
    \gplgaddtomacro\gplbacktext{%
      \csname LTb\endcsname%
      \put(1204,896){\makebox(0,0)[r]{\strut{}$-0.4$}}%
      \put(1204,1848){\makebox(0,0)[r]{\strut{}$-0.2$}}%
      \put(1204,2800){\makebox(0,0)[r]{\strut{}$0$}}%
      \put(1204,3751){\makebox(0,0)[r]{\strut{}$0.2$}}%
      \put(1204,4703){\makebox(0,0)[r]{\strut{}$0.4$}}%
      \put(1372,616){\makebox(0,0){\strut{}$-0.4$}}%
      \put(2703,616){\makebox(0,0){\strut{}$-0.2$}}%
      \put(4034,616){\makebox(0,0){\strut{}$0$}}%
      \put(5364,616){\makebox(0,0){\strut{}$0.2$}}%
      \put(6695,616){\makebox(0,0){\strut{}$0.4$}}%
    }%
    \gplgaddtomacro\gplfronttext{%
      \csname LTb\endcsname%
      \put(224,2799){\rotatebox{-270}{\makebox(0,0){\strut{}parameter $\hat{\theta_2}$}}}%
      \put(4033,196){\makebox(0,0){\strut{}parameter $\hat{\theta_1}$}}%
    }%
    \gplbacktext
    \put(0,0){\includegraphics{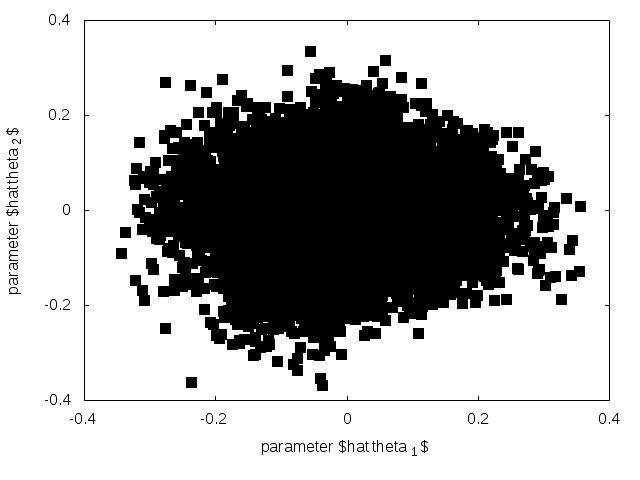}}%
    \gplfronttext
  \end{picture}%
\endgroup

%% file: figures/ExpESS.tex
\begingroup
  \makeatletter
  \providecommand\color[2][]{%
    \GenericError{(gnuplot) \space\space\space\@spaces}{%
      Package color not loaded in conjunction with
      terminal option `colourtext'%
    }{See the gnuplot documentation for explanation.%
    }{Either use 'blacktext' in gnuplot or load the package
      color.sty in LaTeX.}%
    \renewcommand\color[2][]{}%
  }%
  \providecommand\includegraphics[2][]{%
    \GenericError{(gnuplot) \space\space\space\@spaces}{%
      Package graphicx or graphics not loaded%
    }{See the gnuplot documentation for explanation.%
    }{The gnuplot epslatex terminal needs graphicx.sty or graphics.sty.}%
    \renewcommand\includegraphics[2][]{}%
  }%
  \providecommand\rotatebox[2]{#2}%
  \@ifundefined{ifGPcolor}{%
    \newif\ifGPcolor
    \GPcolortrue
  }{}%
  \@ifundefined{ifGPblacktext}{%
    \newif\ifGPblacktext
    \GPblacktexttrue
  }{}%
  \let\gplgaddtomacro\g@addto@macro
  \gdef\gplbacktext{}%
  \gdef\gplfronttext{}%
  \makeatother
  \ifGPblacktext
    \def\colorrgb#1{}%
    \def\colorgray#1{}%
  \else
    \ifGPcolor
      \def\colorrgb#1{\color[rgb]{#1}}%
      \def\colorgray#1{\color[gray]{#1}}%
      \expandafter\def\csname LTw\endcsname{\color{white}}%
      \expandafter\def\csname LTb\endcsname{\color{black}}%
      \expandafter\def\csname LTa\endcsname{\color{black}}%
      \expandafter\def\csname LT0\endcsname{\color[rgb]{1,0,0}}%
      \expandafter\def\csname LT1\endcsname{\color[rgb]{0,1,0}}%
      \expandafter\def\csname LT2\endcsname{\color[rgb]{0,0,1}}%
      \expandafter\def\csname LT3\endcsname{\color[rgb]{1,0,1}}%
      \expandafter\def\csname LT4\endcsname{\color[rgb]{0,1,1}}%
      \expandafter\def\csname LT5\endcsname{\color[rgb]{1,1,0}}%
      \expandafter\def\csname LT6\endcsname{\color[rgb]{0,0,0}}%
      \expandafter\def\csname LT7\endcsname{\color[rgb]{1,0.3,0}}%
      \expandafter\def\csname LT8\endcsname{\color[rgb]{0.5,0.5,0.5}}%
    \else
      \def\colorrgb#1{\color{black}}%
      \def\colorgray#1{\color[gray]{#1}}%
      \expandafter\def\csname LTw\endcsname{\color{white}}%
      \expandafter\def\csname LTb\endcsname{\color{black}}%
      \expandafter\def\csname LTa\endcsname{\color{black}}%
      \expandafter\def\csname LT0\endcsname{\color{black}}%
      \expandafter\def\csname LT1\endcsname{\color{black}}%
      \expandafter\def\csname LT2\endcsname{\color{black}}%
      \expandafter\def\csname LT3\endcsname{\color{black}}%
      \expandafter\def\csname LT4\endcsname{\color{black}}%
      \expandafter\def\csname LT5\endcsname{\color{black}}%
      \expandafter\def\csname LT6\endcsname{\color{black}}%
      \expandafter\def\csname LT7\endcsname{\color{black}}%
      \expandafter\def\csname LT8\endcsname{\color{black}}%
    \fi
  \fi
    \setlength{\unitlength}{0.0500bp}%
    \ifx\gptboxheight\undefined%
      \newlength{\gptboxheight}%
      \newlength{\gptboxwidth}%
      \newsavebox{\gptboxtext}%
    \fi%
    \setlength{\fboxrule}{0.5pt}%
    \setlength{\fboxsep}{1pt}%
\begin{picture}(7200.00,5040.00)%
    \gplgaddtomacro\gplbacktext{%
      \csname LTb\endcsname%
      \put(1372,896){\makebox(0,0)[r]{\strut{}1e+01}}%
      \csname LTb\endcsname%
      \put(1372,2165){\makebox(0,0)[r]{\strut{}1e+02}}%
      \csname LTb\endcsname%
      \put(1372,3434){\makebox(0,0)[r]{\strut{}1e+03}}%
      \csname LTb\endcsname%
      \put(1372,4703){\makebox(0,0)[r]{\strut{}1e+04}}%
      \csname LTb\endcsname%
      \put(1540,616){\makebox(0,0){\strut{}1e+02}}%
      \csname LTb\endcsname%
      \put(4118,616){\makebox(0,0){\strut{}1e+03}}%
      \csname LTb\endcsname%
      \put(6695,616){\makebox(0,0){\strut{}1e+04}}%
    }%
    \gplgaddtomacro\gplfronttext{%
      \csname LTb\endcsname%
      \put(224,2799){\rotatebox{-270}{\makebox(0,0){\strut{}effective sample size}}}%
      \put(4117,196){\makebox(0,0){\strut{}runtime [s]}}%
      \csname LTb\endcsname%
      \put(5456,1379){\makebox(0,0)[r]{\strut{}DRAM}}%
      \csname LTb\endcsname%
      \put(5456,1099){\makebox(0,0)[r]{\strut{}MFMH, independence sampler}}%
    }%
    \gplbacktext
    \put(0,0){\includegraphics{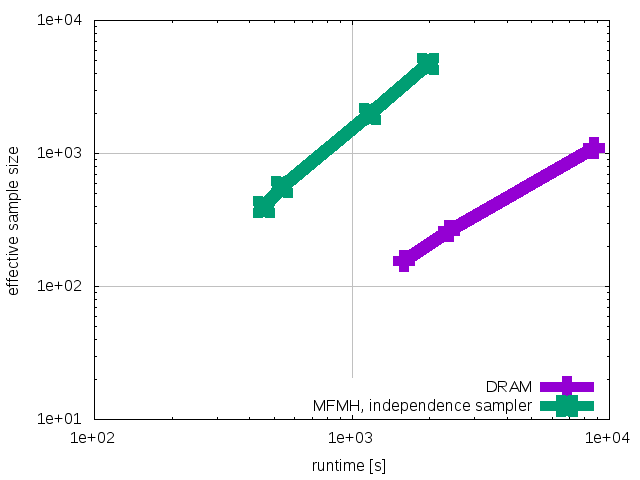}}%
    \gplfronttext
  \end{picture}%
\endgroup

%% file: figures/EulerPosteriorFOM.tex
\begingroup
  \makeatletter
  \providecommand\color[2][]{%
    \GenericError{(gnuplot) \space\space\space\@spaces}{%
      Package color not loaded in conjunction with
      terminal option `colourtext'%
    }{See the gnuplot documentation for explanation.%
    }{Either use 'blacktext' in gnuplot or load the package
      color.sty in LaTeX.}%
    \renewcommand\color[2][]{}%
  }%
  \providecommand\includegraphics[2][]{%
    \GenericError{(gnuplot) \space\space\space\@spaces}{%
      Package graphicx or graphics not loaded%
    }{See the gnuplot documentation for explanation.%
    }{The gnuplot epslatex terminal needs graphicx.sty or graphics.sty.}%
    \renewcommand\includegraphics[2][]{}%
  }%
  \providecommand\rotatebox[2]{#2}%
  \@ifundefined{ifGPcolor}{%
    \newif\ifGPcolor
    \GPcolortrue
  }{}%
  \@ifundefined{ifGPblacktext}{%
    \newif\ifGPblacktext
    \GPblacktexttrue
  }{}%
  \let\gplgaddtomacro\g@addto@macro
  \gdef\gplbacktext{}%
  \gdef\gplfronttext{}%
  \makeatother
  \ifGPblacktext
    \def\colorrgb#1{}%
    \def\colorgray#1{}%
  \else
    \ifGPcolor
      \def\colorrgb#1{\color[rgb]{#1}}%
      \def\colorgray#1{\color[gray]{#1}}%
      \expandafter\def\csname LTw\endcsname{\color{white}}%
      \expandafter\def\csname LTb\endcsname{\color{black}}%
      \expandafter\def\csname LTa\endcsname{\color{black}}%
      \expandafter\def\csname LT0\endcsname{\color[rgb]{1,0,0}}%
      \expandafter\def\csname LT1\endcsname{\color[rgb]{0,1,0}}%
      \expandafter\def\csname LT2\endcsname{\color[rgb]{0,0,1}}%
      \expandafter\def\csname LT3\endcsname{\color[rgb]{1,0,1}}%
      \expandafter\def\csname LT4\endcsname{\color[rgb]{0,1,1}}%
      \expandafter\def\csname LT5\endcsname{\color[rgb]{1,1,0}}%
      \expandafter\def\csname LT6\endcsname{\color[rgb]{0,0,0}}%
      \expandafter\def\csname LT7\endcsname{\color[rgb]{1,0.3,0}}%
      \expandafter\def\csname LT8\endcsname{\color[rgb]{0.5,0.5,0.5}}%
    \else
      \def\colorrgb#1{\color{black}}%
      \def\colorgray#1{\color[gray]{#1}}%
      \expandafter\def\csname LTw\endcsname{\color{white}}%
      \expandafter\def\csname LTb\endcsname{\color{black}}%
      \expandafter\def\csname LTa\endcsname{\color{black}}%
      \expandafter\def\csname LT0\endcsname{\color{black}}%
      \expandafter\def\csname LT1\endcsname{\color{black}}%
      \expandafter\def\csname LT2\endcsname{\color{black}}%
      \expandafter\def\csname LT3\endcsname{\color{black}}%
      \expandafter\def\csname LT4\endcsname{\color{black}}%
      \expandafter\def\csname LT5\endcsname{\color{black}}%
      \expandafter\def\csname LT6\endcsname{\color{black}}%
      \expandafter\def\csname LT7\endcsname{\color{black}}%
      \expandafter\def\csname LT8\endcsname{\color{black}}%
    \fi
  \fi
    \setlength{\unitlength}{0.0500bp}%
    \ifx\gptboxheight\undefined%
      \newlength{\gptboxheight}%
      \newlength{\gptboxwidth}%
      \newsavebox{\gptboxtext}%
    \fi%
    \setlength{\fboxrule}{0.5pt}%
    \setlength{\fboxsep}{1pt}%
\begin{picture}(7200.00,5040.00)%
    \gplgaddtomacro\gplbacktext{%
    }%
    \gplgaddtomacro\gplfronttext{%
      \csname LTb\endcsname%
      \put(1247,4955){\makebox(0,0){\strut{}parameter $\theta_3$}}%
    }%
    \gplgaddtomacro\gplbacktext{%
    }%
    \gplgaddtomacro\gplfronttext{%
    }%
    \gplgaddtomacro\gplbacktext{%
    }%
    \gplgaddtomacro\gplfronttext{%
    }%
    \gplgaddtomacro\gplbacktext{%
      \csname LTb\endcsname%
      \put(1248,3740){\makebox(0,0){\strut{}$\quad$}}%
    }%
    \gplgaddtomacro\gplfronttext{%
    }%
    \gplgaddtomacro\gplbacktext{%
    }%
    \gplgaddtomacro\gplfronttext{%
      \csname LTb\endcsname%
      \put(3599,3460){\makebox(0,0){\strut{}parameter $\theta_2$}}%
    }%
    \gplgaddtomacro\gplbacktext{%
    }%
    \gplgaddtomacro\gplfronttext{%
    }%
    \gplgaddtomacro\gplbacktext{%
      \csname LTb\endcsname%
      \put(1248,2245){\makebox(0,0){\strut{}$\quad$}}%
    }%
    \gplgaddtomacro\gplfronttext{%
    }%
    \gplgaddtomacro\gplbacktext{%
      \csname LTb\endcsname%
      \put(3600,2245){\makebox(0,0){\strut{}$\quad$}}%
    }%
    \gplgaddtomacro\gplfronttext{%
    }%
    \gplgaddtomacro\gplbacktext{%
    }%
    \gplgaddtomacro\gplfronttext{%
      \csname LTb\endcsname%
      \put(5951,1965){\makebox(0,0){\strut{}parameter $\theta_1$}}%
    }%
    \gplbacktext
    \put(0,0){\includegraphics{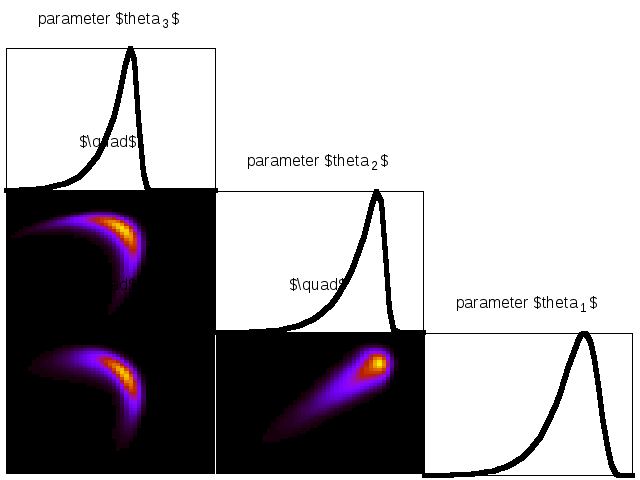}}%
    \gplfronttext
  \end{picture}%
\endgroup

%% file: figures/EulerDRAM.tex
\begingroup
  \makeatletter
  \providecommand\color[2][]{%
    \GenericError{(gnuplot) \space\space\space\@spaces}{%
      Package color not loaded in conjunction with
      terminal option `colourtext'%
    }{See the gnuplot documentation for explanation.%
    }{Either use 'blacktext' in gnuplot or load the package
      color.sty in LaTeX.}%
    \renewcommand\color[2][]{}%
  }%
  \providecommand\includegraphics[2][]{%
    \GenericError{(gnuplot) \space\space\space\@spaces}{%
      Package graphicx or graphics not loaded%
    }{See the gnuplot documentation for explanation.%
    }{The gnuplot epslatex terminal needs graphicx.sty or graphics.sty.}%
    \renewcommand\includegraphics[2][]{}%
  }%
  \providecommand\rotatebox[2]{#2}%
  \@ifundefined{ifGPcolor}{%
    \newif\ifGPcolor
    \GPcolortrue
  }{}%
  \@ifundefined{ifGPblacktext}{%
    \newif\ifGPblacktext
    \GPblacktexttrue
  }{}%
  \let\gplgaddtomacro\g@addto@macro
  \gdef\gplbacktext{}%
  \gdef\gplfronttext{}%
  \makeatother
  \ifGPblacktext
    \def\colorrgb#1{}%
    \def\colorgray#1{}%
  \else
    \ifGPcolor
      \def\colorrgb#1{\color[rgb]{#1}}%
      \def\colorgray#1{\color[gray]{#1}}%
      \expandafter\def\csname LTw\endcsname{\color{white}}%
      \expandafter\def\csname LTb\endcsname{\color{black}}%
      \expandafter\def\csname LTa\endcsname{\color{black}}%
      \expandafter\def\csname LT0\endcsname{\color[rgb]{1,0,0}}%
      \expandafter\def\csname LT1\endcsname{\color[rgb]{0,1,0}}%
      \expandafter\def\csname LT2\endcsname{\color[rgb]{0,0,1}}%
      \expandafter\def\csname LT3\endcsname{\color[rgb]{1,0,1}}%
      \expandafter\def\csname LT4\endcsname{\color[rgb]{0,1,1}}%
      \expandafter\def\csname LT5\endcsname{\color[rgb]{1,1,0}}%
      \expandafter\def\csname LT6\endcsname{\color[rgb]{0,0,0}}%
      \expandafter\def\csname LT7\endcsname{\color[rgb]{1,0.3,0}}%
      \expandafter\def\csname LT8\endcsname{\color[rgb]{0.5,0.5,0.5}}%
    \else
      \def\colorrgb#1{\color{black}}%
      \def\colorgray#1{\color[gray]{#1}}%
      \expandafter\def\csname LTw\endcsname{\color{white}}%
      \expandafter\def\csname LTb\endcsname{\color{black}}%
      \expandafter\def\csname LTa\endcsname{\color{black}}%
      \expandafter\def\csname LT0\endcsname{\color{black}}%
      \expandafter\def\csname LT1\endcsname{\color{black}}%
      \expandafter\def\csname LT2\endcsname{\color{black}}%
      \expandafter\def\csname LT3\endcsname{\color{black}}%
      \expandafter\def\csname LT4\endcsname{\color{black}}%
      \expandafter\def\csname LT5\endcsname{\color{black}}%
      \expandafter\def\csname LT6\endcsname{\color{black}}%
      \expandafter\def\csname LT7\endcsname{\color{black}}%
      \expandafter\def\csname LT8\endcsname{\color{black}}%
    \fi
  \fi
    \setlength{\unitlength}{0.0500bp}%
    \ifx\gptboxheight\undefined%
      \newlength{\gptboxheight}%
      \newlength{\gptboxwidth}%
      \newsavebox{\gptboxtext}%
    \fi%
    \setlength{\fboxrule}{0.5pt}%
    \setlength{\fboxsep}{1pt}%
\begin{picture}(7200.00,5040.00)%
    \gplgaddtomacro\gplbacktext{%
    }%
    \gplgaddtomacro\gplfronttext{%
      \csname LTb\endcsname%
      \put(1247,4955){\makebox(0,0){\strut{}parameter $\theta_3$}}%
    }%
    \gplgaddtomacro\gplbacktext{%
    }%
    \gplgaddtomacro\gplfronttext{%
    }%
    \gplgaddtomacro\gplbacktext{%
    }%
    \gplgaddtomacro\gplfronttext{%
    }%
    \gplgaddtomacro\gplbacktext{%
    }%
    \gplgaddtomacro\gplfronttext{%
      \csname LTb\endcsname%
      \put(1247,3460){\makebox(0,0){\strut{}$\quad$}}%
    }%
    \gplgaddtomacro\gplbacktext{%
    }%
    \gplgaddtomacro\gplfronttext{%
      \csname LTb\endcsname%
      \put(3599,3460){\makebox(0,0){\strut{}parameter $\theta_2$}}%
    }%
    \gplgaddtomacro\gplbacktext{%
    }%
    \gplgaddtomacro\gplfronttext{%
    }%
    \gplgaddtomacro\gplbacktext{%
    }%
    \gplgaddtomacro\gplfronttext{%
      \csname LTb\endcsname%
      \put(1247,1965){\makebox(0,0){\strut{}$\quad$}}%
    }%
    \gplgaddtomacro\gplbacktext{%
    }%
    \gplgaddtomacro\gplfronttext{%
      \csname LTb\endcsname%
      \put(3599,1965){\makebox(0,0){\strut{}$\quad$}}%
    }%
    \gplgaddtomacro\gplbacktext{%
    }%
    \gplgaddtomacro\gplfronttext{%
      \csname LTb\endcsname%
      \put(5951,1965){\makebox(0,0){\strut{}parameter $\theta_1$}}%
    }%
    \gplbacktext
    \put(0,0){\includegraphics{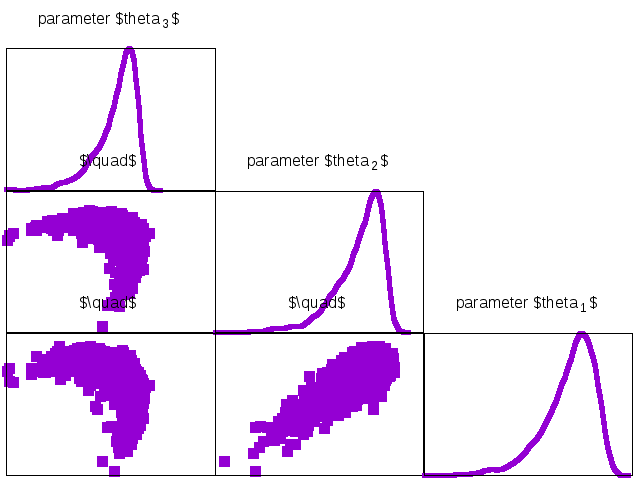}}%
    \gplfronttext
  \end{picture}%
\endgroup

%% file: figures/EulerROMIndep.tex
\begingroup
  \makeatletter
  \providecommand\color[2][]{%
    \GenericError{(gnuplot) \space\space\space\@spaces}{%
      Package color not loaded in conjunction with
      terminal option `colourtext'%
    }{See the gnuplot documentation for explanation.%
    }{Either use 'blacktext' in gnuplot or load the package
      color.sty in LaTeX.}%
    \renewcommand\color[2][]{}%
  }%
  \providecommand\includegraphics[2][]{%
    \GenericError{(gnuplot) \space\space\space\@spaces}{%
      Package graphicx or graphics not loaded%
    }{See the gnuplot documentation for explanation.%
    }{The gnuplot epslatex terminal needs graphicx.sty or graphics.sty.}%
    \renewcommand\includegraphics[2][]{}%
  }%
  \providecommand\rotatebox[2]{#2}%
  \@ifundefined{ifGPcolor}{%
    \newif\ifGPcolor
    \GPcolortrue
  }{}%
  \@ifundefined{ifGPblacktext}{%
    \newif\ifGPblacktext
    \GPblacktexttrue
  }{}%
  \let\gplgaddtomacro\g@addto@macro
  \gdef\gplbacktext{}%
  \gdef\gplfronttext{}%
  \makeatother
  \ifGPblacktext
    \def\colorrgb#1{}%
    \def\colorgray#1{}%
  \else
    \ifGPcolor
      \def\colorrgb#1{\color[rgb]{#1}}%
      \def\colorgray#1{\color[gray]{#1}}%
      \expandafter\def\csname LTw\endcsname{\color{white}}%
      \expandafter\def\csname LTb\endcsname{\color{black}}%
      \expandafter\def\csname LTa\endcsname{\color{black}}%
      \expandafter\def\csname LT0\endcsname{\color[rgb]{1,0,0}}%
      \expandafter\def\csname LT1\endcsname{\color[rgb]{0,1,0}}%
      \expandafter\def\csname LT2\endcsname{\color[rgb]{0,0,1}}%
      \expandafter\def\csname LT3\endcsname{\color[rgb]{1,0,1}}%
      \expandafter\def\csname LT4\endcsname{\color[rgb]{0,1,1}}%
      \expandafter\def\csname LT5\endcsname{\color[rgb]{1,1,0}}%
      \expandafter\def\csname LT6\endcsname{\color[rgb]{0,0,0}}%
      \expandafter\def\csname LT7\endcsname{\color[rgb]{1,0.3,0}}%
      \expandafter\def\csname LT8\endcsname{\color[rgb]{0.5,0.5,0.5}}%
    \else
      \def\colorrgb#1{\color{black}}%
      \def\colorgray#1{\color[gray]{#1}}%
      \expandafter\def\csname LTw\endcsname{\color{white}}%
      \expandafter\def\csname LTb\endcsname{\color{black}}%
      \expandafter\def\csname LTa\endcsname{\color{black}}%
      \expandafter\def\csname LT0\endcsname{\color{black}}%
      \expandafter\def\csname LT1\endcsname{\color{black}}%
      \expandafter\def\csname LT2\endcsname{\color{black}}%
      \expandafter\def\csname LT3\endcsname{\color{black}}%
      \expandafter\def\csname LT4\endcsname{\color{black}}%
      \expandafter\def\csname LT5\endcsname{\color{black}}%
      \expandafter\def\csname LT6\endcsname{\color{black}}%
      \expandafter\def\csname LT7\endcsname{\color{black}}%
      \expandafter\def\csname LT8\endcsname{\color{black}}%
    \fi
  \fi
    \setlength{\unitlength}{0.0500bp}%
    \ifx\gptboxheight\undefined%
      \newlength{\gptboxheight}%
      \newlength{\gptboxwidth}%
      \newsavebox{\gptboxtext}%
    \fi%
    \setlength{\fboxrule}{0.5pt}%
    \setlength{\fboxsep}{1pt}%
\begin{picture}(7200.00,5040.00)%
    \gplgaddtomacro\gplbacktext{%
    }%
    \gplgaddtomacro\gplfronttext{%
      \csname LTb\endcsname%
      \put(1247,4955){\makebox(0,0){\strut{}parameter $\theta_3$}}%
    }%
    \gplgaddtomacro\gplbacktext{%
    }%
    \gplgaddtomacro\gplfronttext{%
    }%
    \gplgaddtomacro\gplbacktext{%
    }%
    \gplgaddtomacro\gplfronttext{%
    }%
    \gplgaddtomacro\gplbacktext{%
    }%
    \gplgaddtomacro\gplfronttext{%
      \csname LTb\endcsname%
      \put(1247,3460){\makebox(0,0){\strut{}$\quad$}}%
    }%
    \gplgaddtomacro\gplbacktext{%
    }%
    \gplgaddtomacro\gplfronttext{%
      \csname LTb\endcsname%
      \put(3599,3460){\makebox(0,0){\strut{}parameter $\theta_2$}}%
    }%
    \gplgaddtomacro\gplbacktext{%
    }%
    \gplgaddtomacro\gplfronttext{%
    }%
    \gplgaddtomacro\gplbacktext{%
    }%
    \gplgaddtomacro\gplfronttext{%
      \csname LTb\endcsname%
      \put(1247,1965){\makebox(0,0){\strut{}$\quad$}}%
    }%
    \gplgaddtomacro\gplbacktext{%
    }%
    \gplgaddtomacro\gplfronttext{%
      \csname LTb\endcsname%
      \put(3599,1965){\makebox(0,0){\strut{}$\quad$}}%
    }%
    \gplgaddtomacro\gplbacktext{%
    }%
    \gplgaddtomacro\gplfronttext{%
      \csname LTb\endcsname%
      \put(5951,1965){\makebox(0,0){\strut{}parameter $\theta_1$}}%
    }%
    \gplbacktext
    \put(0,0){\includegraphics{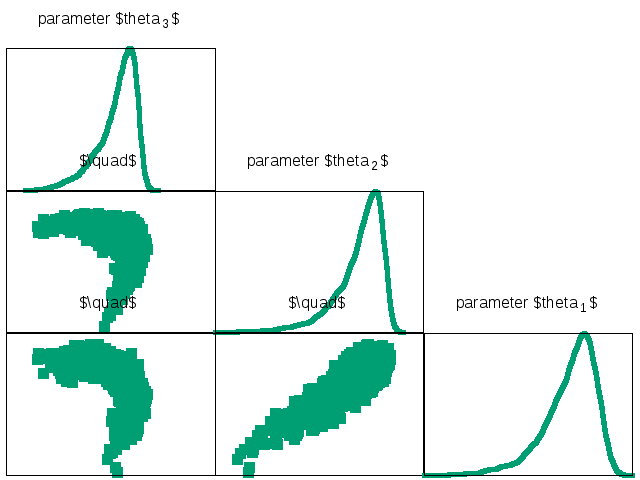}}%
    \gplfronttext
  \end{picture}%
\endgroup

%% file: figures/EulerROMIndepProj.tex
\begingroup
  \makeatletter
  \providecommand\color[2][]{%
    \GenericError{(gnuplot) \space\space\space\@spaces}{%
      Package color not loaded in conjunction with
      terminal option `colourtext'%
    }{See the gnuplot documentation for explanation.%
    }{Either use 'blacktext' in gnuplot or load the package
      color.sty in LaTeX.}%
    \renewcommand\color[2][]{}%
  }%
  \providecommand\includegraphics[2][]{%
    \GenericError{(gnuplot) \space\space\space\@spaces}{%
      Package graphicx or graphics not loaded%
    }{See the gnuplot documentation for explanation.%
    }{The gnuplot epslatex terminal needs graphicx.sty or graphics.sty.}%
    \renewcommand\includegraphics[2][]{}%
  }%
  \providecommand\rotatebox[2]{#2}%
  \@ifundefined{ifGPcolor}{%
    \newif\ifGPcolor
    \GPcolortrue
  }{}%
  \@ifundefined{ifGPblacktext}{%
    \newif\ifGPblacktext
    \GPblacktexttrue
  }{}%
  \let\gplgaddtomacro\g@addto@macro
  \gdef\gplbacktext{}%
  \gdef\gplfronttext{}%
  \makeatother
  \ifGPblacktext
    \def\colorrgb#1{}%
    \def\colorgray#1{}%
  \else
    \ifGPcolor
      \def\colorrgb#1{\color[rgb]{#1}}%
      \def\colorgray#1{\color[gray]{#1}}%
      \expandafter\def\csname LTw\endcsname{\color{white}}%
      \expandafter\def\csname LTb\endcsname{\color{black}}%
      \expandafter\def\csname LTa\endcsname{\color{black}}%
      \expandafter\def\csname LT0\endcsname{\color[rgb]{1,0,0}}%
      \expandafter\def\csname LT1\endcsname{\color[rgb]{0,1,0}}%
      \expandafter\def\csname LT2\endcsname{\color[rgb]{0,0,1}}%
      \expandafter\def\csname LT3\endcsname{\color[rgb]{1,0,1}}%
      \expandafter\def\csname LT4\endcsname{\color[rgb]{0,1,1}}%
      \expandafter\def\csname LT5\endcsname{\color[rgb]{1,1,0}}%
      \expandafter\def\csname LT6\endcsname{\color[rgb]{0,0,0}}%
      \expandafter\def\csname LT7\endcsname{\color[rgb]{1,0.3,0}}%
      \expandafter\def\csname LT8\endcsname{\color[rgb]{0.5,0.5,0.5}}%
    \else
      \def\colorrgb#1{\color{black}}%
      \def\colorgray#1{\color[gray]{#1}}%
      \expandafter\def\csname LTw\endcsname{\color{white}}%
      \expandafter\def\csname LTb\endcsname{\color{black}}%
      \expandafter\def\csname LTa\endcsname{\color{black}}%
      \expandafter\def\csname LT0\endcsname{\color{black}}%
      \expandafter\def\csname LT1\endcsname{\color{black}}%
      \expandafter\def\csname LT2\endcsname{\color{black}}%
      \expandafter\def\csname LT3\endcsname{\color{black}}%
      \expandafter\def\csname LT4\endcsname{\color{black}}%
      \expandafter\def\csname LT5\endcsname{\color{black}}%
      \expandafter\def\csname LT6\endcsname{\color{black}}%
      \expandafter\def\csname LT7\endcsname{\color{black}}%
      \expandafter\def\csname LT8\endcsname{\color{black}}%
    \fi
  \fi
    \setlength{\unitlength}{0.0500bp}%
    \ifx\gptboxheight\undefined%
      \newlength{\gptboxheight}%
      \newlength{\gptboxwidth}%
      \newsavebox{\gptboxtext}%
    \fi%
    \setlength{\fboxrule}{0.5pt}%
    \setlength{\fboxsep}{1pt}%
\begin{picture}(7200.00,5040.00)%
    \gplgaddtomacro\gplbacktext{%
    }%
    \gplgaddtomacro\gplfronttext{%
      \csname LTb\endcsname%
      \put(1247,4955){\makebox(0,0){\strut{}parameter $\theta_3$}}%
    }%
    \gplgaddtomacro\gplbacktext{%
    }%
    \gplgaddtomacro\gplfronttext{%
    }%
    \gplgaddtomacro\gplbacktext{%
    }%
    \gplgaddtomacro\gplfronttext{%
    }%
    \gplgaddtomacro\gplbacktext{%
    }%
    \gplgaddtomacro\gplfronttext{%
      \csname LTb\endcsname%
      \put(1247,3460){\makebox(0,0){\strut{}$\quad$}}%
    }%
    \gplgaddtomacro\gplbacktext{%
    }%
    \gplgaddtomacro\gplfronttext{%
      \csname LTb\endcsname%
      \put(3599,3460){\makebox(0,0){\strut{}parameter $\theta_2$}}%
    }%
    \gplgaddtomacro\gplbacktext{%
    }%
    \gplgaddtomacro\gplfronttext{%
    }%
    \gplgaddtomacro\gplbacktext{%
    }%
    \gplgaddtomacro\gplfronttext{%
      \csname LTb\endcsname%
      \put(1247,1965){\makebox(0,0){\strut{}$\quad$}}%
    }%
    \gplgaddtomacro\gplbacktext{%
    }%
    \gplgaddtomacro\gplfronttext{%
      \csname LTb\endcsname%
      \put(3599,1965){\makebox(0,0){\strut{}$\quad$}}%
    }%
    \gplgaddtomacro\gplbacktext{%
    }%
    \gplgaddtomacro\gplfronttext{%
      \csname LTb\endcsname%
      \put(5951,1965){\makebox(0,0){\strut{}parameter $\theta_1$}}%
    }%
    \gplbacktext
    \put(0,0){\includegraphics{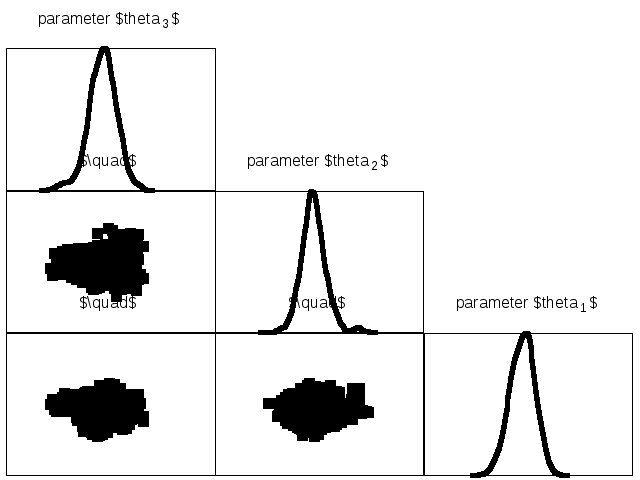}}%
    \gplfronttext
  \end{picture}%
\endgroup

%% file: figures/EulerESS.tex
\begingroup
  \makeatletter
  \providecommand\color[2][]{%
    \GenericError{(gnuplot) \space\space\space\@spaces}{%
      Package color not loaded in conjunction with
      terminal option `colourtext'%
    }{See the gnuplot documentation for explanation.%
    }{Either use 'blacktext' in gnuplot or load the package
      color.sty in LaTeX.}%
    \renewcommand\color[2][]{}%
  }%
  \providecommand\includegraphics[2][]{%
    \GenericError{(gnuplot) \space\space\space\@spaces}{%
      Package graphicx or graphics not loaded%
    }{See the gnuplot documentation for explanation.%
    }{The gnuplot epslatex terminal needs graphicx.sty or graphics.sty.}%
    \renewcommand\includegraphics[2][]{}%
  }%
  \providecommand\rotatebox[2]{#2}%
  \@ifundefined{ifGPcolor}{%
    \newif\ifGPcolor
    \GPcolortrue
  }{}%
  \@ifundefined{ifGPblacktext}{%
    \newif\ifGPblacktext
    \GPblacktexttrue
  }{}%
  \let\gplgaddtomacro\g@addto@macro
  \gdef\gplbacktext{}%
  \gdef\gplfronttext{}%
  \makeatother
  \ifGPblacktext
    \def\colorrgb#1{}%
    \def\colorgray#1{}%
  \else
    \ifGPcolor
      \def\colorrgb#1{\color[rgb]{#1}}%
      \def\colorgray#1{\color[gray]{#1}}%
      \expandafter\def\csname LTw\endcsname{\color{white}}%
      \expandafter\def\csname LTb\endcsname{\color{black}}%
      \expandafter\def\csname LTa\endcsname{\color{black}}%
      \expandafter\def\csname LT0\endcsname{\color[rgb]{1,0,0}}%
      \expandafter\def\csname LT1\endcsname{\color[rgb]{0,1,0}}%
      \expandafter\def\csname LT2\endcsname{\color[rgb]{0,0,1}}%
      \expandafter\def\csname LT3\endcsname{\color[rgb]{1,0,1}}%
      \expandafter\def\csname LT4\endcsname{\color[rgb]{0,1,1}}%
      \expandafter\def\csname LT5\endcsname{\color[rgb]{1,1,0}}%
      \expandafter\def\csname LT6\endcsname{\color[rgb]{0,0,0}}%
      \expandafter\def\csname LT7\endcsname{\color[rgb]{1,0.3,0}}%
      \expandafter\def\csname LT8\endcsname{\color[rgb]{0.5,0.5,0.5}}%
    \else
      \def\colorrgb#1{\color{black}}%
      \def\colorgray#1{\color[gray]{#1}}%
      \expandafter\def\csname LTw\endcsname{\color{white}}%
      \expandafter\def\csname LTb\endcsname{\color{black}}%
      \expandafter\def\csname LTa\endcsname{\color{black}}%
      \expandafter\def\csname LT0\endcsname{\color{black}}%
      \expandafter\def\csname LT1\endcsname{\color{black}}%
      \expandafter\def\csname LT2\endcsname{\color{black}}%
      \expandafter\def\csname LT3\endcsname{\color{black}}%
      \expandafter\def\csname LT4\endcsname{\color{black}}%
      \expandafter\def\csname LT5\endcsname{\color{black}}%
      \expandafter\def\csname LT6\endcsname{\color{black}}%
      \expandafter\def\csname LT7\endcsname{\color{black}}%
      \expandafter\def\csname LT8\endcsname{\color{black}}%
    \fi
  \fi
    \setlength{\unitlength}{0.0500bp}%
    \ifx\gptboxheight\undefined%
      \newlength{\gptboxheight}%
      \newlength{\gptboxwidth}%
      \newsavebox{\gptboxtext}%
    \fi%
    \setlength{\fboxrule}{0.5pt}%
    \setlength{\fboxsep}{1pt}%
\begin{picture}(7200.00,5040.00)%
    \gplgaddtomacro\gplbacktext{%
      \csname LTb\endcsname%
      \put(1372,1336){\makebox(0,0)[r]{\strut{}1e+04}}%
      \csname LTb\endcsname%
      \put(1372,2800){\makebox(0,0)[r]{\strut{}1e+05}}%
      \csname LTb\endcsname%
      \put(1372,4263){\makebox(0,0)[r]{\strut{}1e+06}}%
      \csname LTb\endcsname%
      \put(1540,616){\makebox(0,0){\strut{}1e+03}}%
      \csname LTb\endcsname%
      \put(4118,616){\makebox(0,0){\strut{}1e+04}}%
      \csname LTb\endcsname%
      \put(6695,616){\makebox(0,0){\strut{}1e+05}}%
    }%
    \gplgaddtomacro\gplfronttext{%
      \csname LTb\endcsname%
      \put(224,2799){\rotatebox{-270}{\makebox(0,0){\strut{}effective sample size}}}%
      \put(4117,196){\makebox(0,0){\strut{}runtime [s]}}%
      \csname LTb\endcsname%
      \put(5456,1659){\makebox(0,0)[r]{\strut{}DRAM}}%
      \csname LTb\endcsname%
      \put(5456,1379){\makebox(0,0)[r]{\strut{}MFMH, local random walk}}%
      \csname LTb\endcsname%
      \put(5456,1099){\makebox(0,0)[r]{\strut{}MFMH, independence sampler}}%
    }%
    \gplbacktext
    \put(0,0){\includegraphics{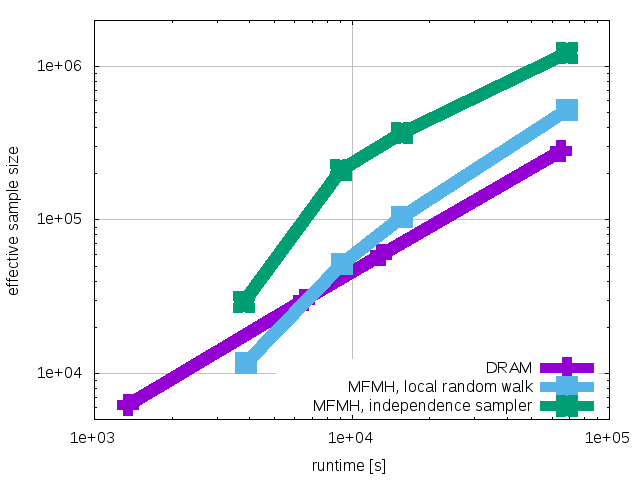}}%
    \gplfronttext
  \end{picture}%
\endgroup